\documentclass[12pt]{amsart}

\usepackage{pdfsync}
\usepackage{hyperref}
 \usepackage[normalem]{ulem}

\usepackage{color}

 \let\cal\mathcal

\usepackage{graphicx}
\usepackage[curve]{xypic}
\usepackage[usenames,dvipsnames,x11names]{xcolor}
\usepackage{tikz}

\newtheorem{theorem}{Theorem}[section]
\newtheorem{proposition}[theorem]{Proposition}
\newtheorem{corollary}[theorem]{Corollary}
\newtheorem{lemma}[theorem]{Lemma}
\newtheorem*{lemma*}{Lemma}

\theoremstyle{definition}
\newtheorem*{definition*}{Definition}
\newtheorem{definition}[theorem]{Definition}

\newtheorem{remark}[theorem]{Remark}
\newtheorem*{remark*}{Remark}

\newtheorem{rien}[theorem]{}
\newtheorem*{theorem*}{Theorem}
\newtheorem*{proposition*}{Proposition}

\usepackage{graphicx}
\begin{document}
\title{ Carrousel in family and non-isolated hypersurface singularities  in $\Bbb C^3$  }
\author{Fran\c{c}oise Michel }
 \author{Anne  Pichon}
\maketitle

\begin{abstract}
 We study the boundary $L_t$ of the Milnor fiber for the
reduced holomorphic germs $f:(\Bbb C^3,0) \rightarrow (\Bbb C,0)$ having a non-isolated singularity at $0$.  We prove that $L_t$ is a graph manifold by using a new technique of carrousels  depending on one  parameter. Our results enable us to compare the topology of $L_t$ and of the link of the normalization of $f^{-1}(0)$.
\end{abstract}
   
\section { Introduction} 
 
  We denote by $ \Bbb B^{2n}_r$  the $2n$-ball
with radius $r>0$ centered at the origin of ${\Bbb C}^n$ and by
$\Bbb S^{2n-1}_r$  the boundary of $\Bbb B^{2n}_r$.
 
Let   $f : { (\Bbb C}^3,0) \longrightarrow {(\Bbb C},0)$ be a reduced  holomorphic germ. The singularity of $f$ at $0$ is allowed to be non-isolated.  We consider  the  three underlying topological objects: 
 
 \begin{itemize}
\item[$\bullet$] The link $L_0 = f^{-1}(0) \cap \Bbb S^5_{\epsilon}$ of the   surface $ f^{-1}(0)$ at $0$, whose homeomorphism class does not depend on $\epsilon$ when $\epsilon >0$ is  sufficiently small (\cite{Mi}, \cite{B-V}).
\item[$\bullet$] The boundary  $L_t =  f^{-1}(t) \cap \Bbb S^5_{\epsilon}$ of the Milnor fiber of $f$, where $0 < |t|   << \epsilon$, whose diffeomorphism class does not depend on $t$ when $|t|$ is sufficiently small (\cite{Mi}, \cite{H-L}). 
\item[$\bullet$] The link $\overline{L_0}$ of the normalization of    the surface $F_0 =  f^{-1}(0) \cap \Bbb B^6_{\epsilon}$ at $0$, which can be defined up to diffeomorphism  by $\overline{L_0} = n^{-1}(L_0)$, where $n : \overline{F_0} \rightarrow F_0$ denotes the
normalization morphism of $F_0$ (\cite{D}).
\end{itemize}

When the origin is an isolated singular point, $L_0$, $L_t$ and $\overline{L_0}$ are $3$-dimensional  differentiable
manifolds, each of them being    diffeomorphic  to the others. 

 In  this paper,  we assume that the singular locus  $\Sigma (f)$ of $f$ is   $1$-dimensional. Then only $L_t$ and $\overline{L_0}$ are differentiable manifolds.

Resolution    theory   implies that $\overline{L_0}$ is  a graph manifold in the sense of Waldhausen (\cite{W}),  or equivalently a plumbed manifold  (\cite{HNK}, \cite{N}). More precisely, the plumbing graph of $\overline{L_0}$ is given, in its normal form,   as the dual graph of  a good minimal resolution of the normal surface singularity $ \overline{F_0}$. 

We will not recall here the notions of Seifert, graph and plumbed manifolds.  For a quick survey adapted to our situation, see  e.g. \cite[Section 3]{M-P-W-2}.

{In \cite{M-P} and its erratum, we state  that for a germ $f : { (\Bbb C}^3,0) \longrightarrow {(\Bbb C},0)$,  the boundary $L_t$ of the Milnor fiber is also a graph manifold whose Seifert pieces have oriented basis.   The main aim of this paper is to give a detailed proof of this result. 

We first describe  the manifold $L_t$ using the following strategy (Section 2): by hypothesis the singular locus $\Sigma (f)$ of $f$ is a curve. Let $K_0 = L_0 \cap \Sigma (f)$ be   its link   and let $\overline{K_0} = n^{-1}
(K_0)$ be the pull-back of $K_0$ in $\overline{L_0}$.  Let $\overline  \Sigma (f) = n^{-1} (\Sigma (f))$  be the pull-back of   $\Sigma(f)$ by $n$. A good
resolution of  the pair $(  \overline{F_0}, \overline  \Sigma (f) ) $ provides a    decomposition
for $\overline{L_0}$ as a union of Seifert manifolds such that
$\overline{K_0}$ is a union of  Seifert  fibers. Let $\overline  M_0$ be
a tubular neighborhood of
$\overline{K_0}$ in $\overline{L_0}$. The closure $\overline  N_0$ of
$(\overline{L_0}  \setminus    \overline{M_0})$ is an irreducible  graph manifold with boundary that we  called the {\it trunk} of $ \overline{L_0}$.

On the other hand, we define (Definition \ref{def:vanishingZ})   a 
submanifold $M_t$ of $L_t$ called the
{\it vanishing zone} around $K_0$. In the  literature, there exists a vanishing homology defined by D. Siersma in \cite{S}, not to be confused with the vanishing zone introduced in  \cite{M-P}. Theorem \ref{main1} (3), Theorem \ref{main} and  Proposition \ref{irred} can be summarized in the following theorem: 

\vskip.1in\noindent 
{\bf Theorem.} 
\begin{enumerate}
\item The closure $N_t$ of $L_t \setminus  
M_t$ is orientation preserving  diffeomorphic  to the trunk $\overline{N_0}$.
\item The
manifold $M_t$ is an irreducible    graph  manifold whose Seifert pieces have oriented basis.
\end{enumerate}
    \vskip.1in   
    
    As a consequence, we obtain (Section \ref{vertical}) that the vertical monodromy introduced by D. Siersma  in \cite{S} is a quasi-finite diffeomorphism. 
  \vskip.1in       
The main aim in the study of the topological aspects of singularities consists of  describing the analytical properties of a singularity which can be characterized through some  underlying topological  objects. One of the most important results in this direction is the following famous theorem of Mumford, which gives a topological characterization of a smooth point on a normal surface: 

\begin{theorem*} (\cite{Mu}) Let (X,0) be the germ of a normal complex surface. If the link  $ L_0$ of $(X,0)$ has the homotopy type of the $3$-sphere, then $0$ is a smooth point of $X$. 
\end{theorem*}

Our description of $L_t$ enables one to compare the homeomorphism classes of $L_t$ and $\overline{L_0}$. This leads to the following  topological characterization of isolated singularities for  holomorphic reduced germs  $f : { (\Bbb C}^3,0) \longrightarrow {(\Bbb C},0)$. 

 \vskip.1in\noindent
{\bf Theorem   \ref{maintop}.}   Let $f:({\Bbb C}^3,0) \longrightarrow ({\Bbb C},0)$ be a reduced holomorphic
germ. Unless $f$ is irreducible and $L_t$ is a lens space, the following assertions are equivalent:
\begin{itemize}
\item[(i)] $f$ is  either smooth or has an isolated singularity at $0$.
\item[(ii)] The boundary $L_t, \ t \not=0$, of the Milnor fibre of $f$ is homeomorphic to the link $\overline{L_0}$ of the
normalization of $f^{-1}(0)$.
\end{itemize}
 
    \vskip.1in   
 When $f$ is reducible, $f$ has a non-isolated singularity and the number of connected components of $\overline{L_0}$ equals, by definition,  the number of irreducible  factors of $f$. In this case, the result immediately follows from the fact that $L_t$ is always connected (Corollary \ref{connected}).
 
The case when $f$ is irreducible and $L_t$ is a lens space  remains open, but we  show that this  concerns a very special  family of singularities: 

 \vskip.1in\noindent
{\bf Proposition   \ref{lens}.}  Let $f:({\Bbb C}^3,0) \longrightarrow ({\Bbb C},0)$ be a reduced holomorphic
germ such that $f$ is irreducible and $L_t$ is  a lens space. Then 
\begin{enumerate}
\item The trunk $\overline{N_0}$ is a solid torus,  $\overline{L_0}$ is a lens space, $\overline{\Sigma}(f) $ is an irreducible  curve germ  and  the minimal  resolution graph of the   pair $(\overline{F_0},\overline{\Sigma}(f))$ is a bamboo with an arrow at one of its extremities,  
 \item $M_t$ is connected with a connected boundary.
\end{enumerate}

\vskip.1in
In \cite{M-P-W-2}, for the germs with equations  $z^m -
g(x,y) = 0$ where $m \geq 2$ and $g(x,y)=0$  is a non-reduced plane curve germ, we proved that $L_t$ is never homeomorphic to $\overline{L_0}$ even if  $L_t$ is a lens space, and that  latter case arises if and only if $m=2$ and $g$ has the analytic type of $xy^l$. 

   In fact, in most of the known cases, $L_t$ is not orientation preserving homeomorphic  to the link $L_X$ of any  complex normal surface singularity $(X,p)$. This happens for various reasons.  In \cite{M-P}, we show that the germ $f(x,y,z)=xy$ has $L_t\cong \Bbb S^2 \times \Bbb S^1$, which is not an irreducible $3$-manifold.  In \cite{M-P-W-1}, we show that for the germs $z^m-x^ky^l=0$ such that $L_t$ is not a lens space, ({\it i.e.,} $(m,k)\neq (2,1)$), the boundary $L_t$  is not an $L_X$ as the intersection form associated to its normalized plumbing graph is never negative definite.  In \cite{M-P-W-2},  we show that for the germ $z^2 - (x^2 -y^3)y^l$, $l$ odd, the boundary $L_t$ is homeomorphic to the boundary of an $L_X$,
   but with the reversed orientation. In \cite{N-S}, A. N\'emethi and A.  Szil\'ard describe the boundary of the Milnor fiber for other  families of  examples. In particular, they obtain some examples in which some edges of the normalized plumbing graph of $L_t$ have a sign $\epsilon = -1$, which never happens for   an $L_X$.

 \vskip.1in

  In order to prove Theorem \ref{main}, we develop two independent techniques which may have their own interest.

  The first one (in Section 3) leads to the following resolution in family theorem for a family of plane curve germs parametrized by a punctured disc (hence the set of parameters is not contractible).  Let us  consider a reduced holomorphic  germ $f : (\Bbb C^3,0) \rightarrow (\Bbb C,0)$. Let  $\alpha>0 $  be such that the polydisc $B(\alpha) :=  \Bbb B^2_{\alpha}\times \Bbb B^2_{\beta} \times \Bbb B^2_{\gamma} $  is {\it a Milnor polydisc for $f$} (Definition \ref{def:polydisc}).    Let $\sigma$ be an irreducible component of $\Sigma(f)$ and let us consider   $ \sigma ^{*}=\sigma \cap ( int(B(\alpha) \setminus   \{x=0\}))$. 
  
 \vskip.1in\noindent
{\bf Theorem \ref{th:resolution}.}  
 There exists a sufficiently small  $\alpha >0$, an  open analytic manifold $ V,$ neighborhood of  $ \sigma ^{*}$ in $ int(B(\alpha) \setminus   \{x=0\}),$  and a composition of a finite number of  blow-ups  along  punctured discs $\pi :  \tilde V \ \rightarrow   V$   starting with the blow-up of $V$ along $\sigma^*$  such that:
\begin{enumerate}
\item \label{res1} $\pi$ restricted on $\tilde{V} \setminus \pi ^{-1} (\sigma ^{*})$ is an isomorphism.
\item  \label{res2} $E=\pi ^{-1} (\sigma ^{*})$ is an analytic normal crossing divisor.
\item \label{res3} For each $p=(x,y,z)\in \sigma ^*$,   $\pi$ restricted on    $ \pi ^{-1}(V \cap ( \{ x \} \times \Bbb C^2) ) $  is an embedded resolution of  the plane  curve germ, at $p$,  $ V \cap ( \{ x \} \times \Bbb C^2) \cap \{ f=0\}$.
\end{enumerate}

The second technique is a carrousel in family parametrized by $x$ varying along a circle.  The ``carrousel"   was introduced by D.T. L\^e  in \cite[p.163]{Le-2} and \cite{Le} to obtain a geometric proof of the monodromy theorem. In \ref{reduction}, we show that it is sufficient to prove that   the  vanishing zone of  $f$  along $\sigma$, defined in  \ref{def:vanishZ},  is a graph manifold when $\sigma$ is the $x-$axis. In Subsections 4.1 and 4.2,  we consider the  family  of map germs $f_{a}(y,z)=f(a,y,z)$ at $(a,0,0)\in  \{a\} \times \Bbb C^2$, obtained  as  hyperplane  sections by $\{x=a\}.$ The  carrousel construction,  for  the germ $f_{a}$ and  the direction $y$,  is based on the study of the discriminant  $\Delta ^{(a)}$ which  is the set of singular values of the morphism $\Psi^{(a)}(y,z)=(a,y,f_{a}(y,z))$. A simultaneous  Puiseux parametrization  allows us to perform a carrousel construction for  the  family of discriminant curves $\Delta ^{(x)}$ where $x$  is varying along $\Bbb S_{\alpha}^1.$    It  is what we mean  by ``carrousel in family".  In Subsection   4.3,  we explain how these   constructions induce a structure of graph manifold on the vanishing zone.

  \vskip.1in

 The fact that $L_t$ is a graph manifold was also stated and proved  by A. Nem\'ethi and A. Szilard in \cite{N-S} with a radically different approach. It also appears in  the preprint  \cite{FB-MN} by J.  Fernandez de Bobadilla and A. Menegon-Neto. 

We acknowledge  N. A'Campo, J.  Fernandez de Bobadilla and D. Massey to have independently  pointed out to us  the gap in \cite{M-P}. 

We are very grateful to the referee for valuable comments which enabled us to improve the redaction of the paper. }

\section{ The  trunk and the vanishing zone} 

In this section, we define the trunk and the vanishing zone of $L_t$. As a preliminary, we start in Subsections \ref{axis1} and  \ref{axis2}  by performing generic choices of the coordinates axis in $\Bbb C^3$. 

\begin{rien} \label{axis1}  {The Weierstrass preparation theorem  implies that we can suppose that  $f$ is a unitary polynomial  in  $ \Bbb C \{x,y\} [z]$}. Then, the intersection 

$$\Gamma _0 := \{ f=0 \} \cap  \{ {{\partial f}\over {\partial z}}  =0 \} $$ 
is  a curve which contains $\Sigma (f)$. 

{The  choice of $f$ in  $ \Bbb C \{x,y\} [z]$ is  convenient. In particular,  it   makes obvious  the inclusion  $(1_t)$ in  Subsection 2.10  and it  simplifies  the proof of Proposition 2.12. }

\vskip.1in
\noindent
{\bf Claim.} For a generic choice  of the $x$-axis,    
$\{  {{\partial f}\over {\partial z}} =0 \} \cap \{ {{\partial f}\over {\partial y}}  =0 \}  $  does not meet the boundary of the Milnor fiber and:

$$ \Sigma (f)=  \Gamma _0  \cap \{  {{\partial f}\over {\partial y}} =0  \}. $$

\begin{proof}
 D.T. L\^e and B. Teissier  (for example see (2.2.2) in \cite {Le} or IV.1.3.2 p.420 in \cite{T})  have proved that,  for a generic  choice of  the $x$-axis, 
$$ \{ {{\partial f}\over {\partial z}} =0 \} \cap \{ {{\partial f}\over {\partial y}} =0 \} = \Sigma(f) \cup \Gamma _{(x,f)},$$  
where the irreducible  components of $ \Gamma _{(x,f)} $ are $1$-dimensional and not included in $\{ f=0\} $. (They  have called $ \Gamma _{(x,f)} $ the {\it  polar curve of $f$ for the direction $x$}.)  So the boundary of the Milnor fiber $f^{-1}(t)\cap {\Bbb B}_{\epsilon}^6$ does not meet $ \Gamma _{(x,f)}$ for $t$ sufficiently small  (but its interior does). \end{proof}

\end{rien} 
\begin{rien} \label{axis2} Let $P : {\Bbb C}^3 \longrightarrow {\Bbb C}^2$ be the map  defined by $$P (x,y,z)= (x,y).$$

The curve  $\Delta _0:=P(\Gamma _0)$  is the {\it discriminant curve}.   After performing a linear change of coordinates in $\Bbb C^2$ if necessary, we can assume that the $x$-axis  is transverse to $\Delta _0$ at the origin and that, in $\Bbb C^3$, the hyperplanes $X_a= \{ x=a \} $ meet $\Gamma_0$ transversely  around the origin.
\end{rien}

\begin{rien} \label{polydisc} For technical reasons, we replace in this paper the standard Milnor  ball $\Bbb B^6_{\epsilon}$ by a polydisc
$$B(\alpha) := \Bbb B^2_{\alpha} \times  \Bbb B^2_{\beta} \times
\Bbb B^2_{\gamma} = \{(x,y,z) \in \Bbb B^6_{\epsilon} , \ \  |x|\leq \alpha,
|y| \leq \beta, |z| \leq \gamma\}$$
where
$0 < \alpha  < \beta  <\gamma  <\epsilon \slash 3 $.
\end{rien}

\begin{definition} \label{def:polydisc} The polydisc $B(\alpha)$ is {\it a
Milnor polydisc for} $f$ if  for each $\alpha '$ with $0 < \alpha ' \leq
\alpha$, 

 \begin{enumerate}
\item  the pair  $(B(\alpha ') ,f^{-1}(0) \cap B(\alpha '))$
is diffeomorphic to the pair
$(\Bbb B^6_{\epsilon} ,f^{-1}(0) \cap \Bbb B^6_{\epsilon})$, 

\item  there exists $\eta$ with $0 < \eta << \alpha '$ such
that:
 \begin{enumerate}
 \item the restriction of $f$ to
$W(\alpha ', \eta) = B(\alpha ') \cap f^{-1} (\Bbb B^2_{\eta} \setminus
\{ 0 \})$ is a locally trivial differentiable fibration over 
$\Bbb B^2_{\eta} \setminus \{ 0 \}$,

\item the isomorphism class of this fibration does not depend   on $\alpha '$ and $\eta$. 
\end{enumerate}
\end{enumerate}
\end{definition}
 \vskip.1in
  We now show that a Milnor polydisc for $f$ exists. We can choose 
$0 < \alpha  < \beta  <\gamma  <\epsilon \slash 3 $  such that the two following inclusions hold:

\begin{itemize}  
\item[($1_0$)]  ${ \partial B(\alpha)} \cap  f^{-1}(0)  \subset  \{ |z| <  \gamma \}$,  
\item[($2_0$)]  $ (\Gamma_0 \cap { \partial B(\alpha)}) \subset    \Bbb S^1_{\alpha} \times int( \Bbb B^2_{\beta}) \times int( \Bbb B^2_{\gamma})$.
\end{itemize}

 {It is then not hard to see that the generic choice of coordinates of Subsections \ref{axis1} and \ref{axis2} and the above choice of  $\alpha ,\beta, \gamma$ imply that  the polydisc $B(\alpha )$ is a Milnor polydisc for $f$ (see e.g.,  Section 1 of  \cite{Le}).}
  \vskip.1in
 In the sequel, we will then replace the objects defined in the introduction by the following: 
 
\begin{itemize}
\item[$\bullet$] For $0 \leq |t| \leq \eta$, 
$$F_t =  f^{-1}(t) \cap B(\alpha) \hbox{ and } L_t = F_t \cap  { \partial B(\alpha)},$$
\item[$\bullet$] $\overline{L_0} = n^{-1}(L_0)$, where $n : \overline{F_0} \rightarrow F_0$ denotes the normalization of $F_0$,
\item[$\bullet$] $K_0 = \Sigma(f) \cap L_0 \hbox{ and } \overline{K_0} = n^{-1}(K_0)$.
\end{itemize}

  \begin{remark} The restriction  of $P$  on $L_0 =  {\partial B(\alpha)} \cap  f^{-1}(0) $  is a ramified cover  whose  ramification locus is the algebraic link  $\Delta _0 \cap \partial(\Bbb B^2_{\alpha} \times {\Bbb B}^2_{\beta})$  and whose  generic    degree  is the  degree  of $f$ in $z$.
   \end{remark}
   
    The above construction implies the following proposition.
   
    \begin{proposition} \label{tube} {Let $V$ be   a sufficiently small tubular neighborhood  of $\Delta _0 \cap \partial(\Bbb B^2_{\alpha} \times {\Bbb B}^2_{\beta})$ in $\partial(\Bbb B^2_{\alpha} \times {\Bbb B}^2_{\beta})$.} The two following conditions  hold:

\begin{enumerate} \item $ V  \subset   \Bbb S^1_{\alpha} \times int( \Bbb B^2_{\beta})$. 
\item Let $M_0$ be the union of the connected components of     $L_0 \cap P^{-1}(V)$ which contain components of the link $K_0$. Then $ \overline{M_0} = n^{-1}(M_0)$ is a tubular neighborhood of $ \overline{K_0} $ in   $\overline{L_0}$.  \hfill $\Box$
\end{enumerate} 
   \end{proposition}
 
    \begin{definition} The {\it trunk} of    $L_0$  is the closure  $N_0$   of $L_0 \setminus M_0$ in $L_0$. The  {\it trunk} of  $ \overline{L_0}$  is the closure  $ \overline{N_0}$   of $ \overline{L_0} \setminus  \overline{M_0}$ in $\overline{L_0}$.    
   \end{definition}

\begin{figure}[ht]
%\centering
\begin{tikzpicture} 

\begin{scope}[xshift=2cm]
%NORMALISATION
\begin{scope}[yshift=10cm]

 \draw[line width=1pt, lightgray  ] 
(5,3.5)  .. controls (5.5,3.5)  and  (5.1,2.8) .. (5.1,2.6) ;  
 \draw[line width=1pt, lightgray  ] 
(5.8,1.2) .. controls (5.8,1.7) and (5.1,2.3)  .. (5.1,2.6) ;
%vanish
 \draw[line width=2pt] 
(5.19,2.95)  .. controls (5.08,2.65)  and  (5.08,2.6) .. (5.23,2.3) ;

 \draw[line width=1pt, lightgray   ] 
(4.7,2.3) .. controls (4.7,2.6) and (4.5,3.5)  .. (5,3.5) ;

\draw[line width=1pt, lightgray  ] 
(4.7,2.3) .. controls (4.7,2) and (4.2,1.6)  .. (4.2,1.1) ;

%vanish
 \draw[line width=2pt] 
(4.67,2.7)  .. controls (4.73,2.5)  and  (4.75,2.2) .. (4.52,1.9) ;

4.7,2.3

\draw[line width=1pt, lightgray  ] 
(5,0) .. controls (5.7,0.1) and (5.8,0.7)  .. (5.8,1.2) ;

 \draw[line width=1pt , lightgray ]
(4.2,1.1) .. controls (4.2,1) and (4,0)  .. (5,0) ;

% Gamma DANS NORMALISATION

%\draw[line width=0.5pt] 
%(1,1.5) .. controls (2,1.8) and (3.3,1.7)  .. (4.3,1.5) ;

\draw[line width=0.5pt] 
(1,1.5) .. controls (2,2.5) and (4,3.4)  .. (5,3.51) ;

\draw[line width=0.5pt] 
(1,1.5) .. controls (2,2) and (4,2.3)  .. (4.7,2.3) ;

\draw[line width=0.5pt] 
(4.65,2.62) .. controls (4.9,2.65) and (5,2.65)  .. (5.1,2.64) ;

\draw[dotted, line width=0.6pt] 
(1,1.5) .. controls (2,2.1) and (4,2.64)  .. (4.65,2.62) ;

\draw[line width=0.5pt] 
(1,1.5) .. controls (2,1.5) and (4.3,0.3)  .. (4.5,0.1) ;

\draw[fill=white] (4.7,2.3)circle(1.5pt);
\draw[fill=white] (5.11,2.64)circle(1.5pt);

 \node at(0.8,1.5)[ ]{$0$};

  \node at(6.1,1.2)[lightgray  ]{$\overline{N_0}$};
  
   \node at(6,2.3)[right]{$\overline{M_0}$};
 \draw[line width=0.3pt, lightgray]   (6,2.3)--(5.25,2.85);
  \draw[line width=0.3pt, lightgray]   (6,2.25)--(4.7,2);
       
       \node at(2.6,3.5)[left]{$\overline{K_0}$};
        \draw[line width=0.3pt, lightgray]   (2.5,3.35)--(4.6,2.4);
         \draw[line width=0.3pt, lightgray]   (2.5,3.4)--(5,2.7);
        
  \node at(3.5,1.2)[left]{$\overline{F_0}$};
  
  \end{scope}

%FLECHE n
 \begin{scope}[yshift=4.5cm]
  \node at(3.1,4.8)[right]{$n$};
  \draw[line width=0.5pt]   (3,5.2)--(3,4.5);
   \draw[line width=0.5pt]   (3.1,4.7)--(3,4.5);
     \draw[line width=0.5pt]   (2.9,4.7)--(3,4.5);
     \end{scope}

 %SURFACE F0  
  \begin{scope}[yshift=5cm]
 %\node at(3.5,1.2)[left]{$ {F_0}$};

 \draw[line width=1pt, lightgray  ]
(4.2,1.1) .. controls (4.2,1) and (4,0)  .. (5,0) ;

\draw[line width=1pt, lightgray ] 
(5,2.5) .. controls (4.8,2.2) and (4.2,1.6)  .. (4.2,1.1) ;

 \draw[line width=1pt, lightgray    ] 
(5,0) .. controls (6,0.1) and (6,1.8)  .. (5,2.5) ;

 \draw[line width=1pt, lightgray   ] 
(5,2.5) .. controls (4.5,2.9) and (4.6,3.5)  .. (5,3.5) ;
%vanish
\draw[line width=2pt] 
(5,2.5) .. controls (4.9,2.6) and (4.8,2.7)  .. (4.75,2.8) ;
\draw[line width=2pt] 
(5,2.5) .. controls (5.1,2.42) and (5.05,2.47)  .. (5.29,2.23) ;
 
 \draw[line width=1pt, lightgray    ] 
(5,3.5)  .. controls (5.2,3.5)  and  (5.4,3.2) .. (5,2.5) ;

%vanish
\draw[line width=2pt] (5.17,2.85)  .. controls (5.05,2.6)  and  (5.05,2.6) .. (5,2.5) ;
\draw[line width=2pt] (5,2.5) .. controls (4.95,2.4) and (4.85,2.3)  .. (4.72,2.15) ;

% Gamma 
\draw[line width=0.5pt] 
(1,1.5) .. controls (2,2.5) and (4,3.4)  .. (5,3.5) ;

\draw[line width=0.5pt] 
(1,1.5) .. controls (2,2) and (4,2.5)  .. (5,2.5) ;

\draw[line width=0.5pt] 
(1,1.5) .. controls (2,1.5) and (4.3,0.3)  .. (4.5,0.1) ;

  \draw[fill=white] (5,2.5)circle(1.5pt);

 \node at(0.8,1.5)[ ]{$0$};
 
  \node at(6.1,1.2)[lightgray  ]{$N_0$};
  
   \node at(6,2.3)[right]{$M_0$};
 \draw[line width=0.3pt, lightgray]   (6,2.3)--(5.13,2.65);

       \node at(2.6,3.5)[left]{$K_0$};
        \draw[line width=0.3pt, lightgray]   (2.5,3.4)--(4.85,2.54);
  
   \node at(0.8,1.5)[ ]{$0$};
  \node at(3.6,1.8)[ ]{$F_0$};
   \node at(2,2.5)[ ]{$\gamma_1$};
    \node at(2.8,2.3)[ ]{$\gamma_2$};
    % \node at(2.8,1.85)[ ]{$\gamma_3$};
      \node at(3,1.1)[ ]{$\gamma_3$};
     % \node at(5.9,0.5)[ ]{$L_0$};
        %\node at(6.5,1.5)[right ]{$\xrightarrow{P}$};

   \node at(8,2.5)[right ]{$\Sigma(f)=\gamma_2 $};
    \node at(8,2)[ right]{$\Gamma_0= \bigcup_{i=1}^3 \gamma_i$  };
     \node at(8,1.5)[right ]{$L_0= M_0 \cup N_0 $ };

\end{scope}
 
%FLECHE P
\begin{scope}[yshift=-0.5cm]
  \node at(3.1,4.8)[right]{$P_{\mid F_0}$};
  \draw[line width=0.5pt]   (3,5.2)--(3,4.5);
   \draw[line width=0.5pt]   (3.1,4.7)--(3,4.5);
     \draw[line width=0.5pt]   (2.9,4.7)--(3,4.5);
     \end{scope}
     
%DELTA_0
 \begin{scope}[xshift=2.5cm]
 
 \node at(0.8,1.5)[ ]{$0$};

\draw[line width=0.5pt] 
(1,1.5) .. controls (2,2.5) and (4,3.4)  .. (5,3.5) ;

\draw[line width=0.5pt] 
(1,1.5) .. controls (2,2) and (4,2.5)  .. (5,2.5) ;

 \draw[line width=0.5pt] 
(1,1.5) .. controls (2,1.5) and (4.3,0.3)  .. (4.5,0.1) ;

\draw[line width=0.3pt] 
(-0.5,3.5)--(-0.5,-0.5)--(2.5,-0.5)--(2.5,3.5)--(-0.5,3.5);

\draw[line width=2pt] 
 (2.5, 2.4)--(2.5,2.8);
 
 \draw[line width=2pt] 
 (2.5, 2.2)--(2.5,1.9);
 
   \draw[line width=2pt] 
 (2.5, 1.3)--(2.5,0.8);

   \node at(3.4,3.25)[ ]{$\delta_1$};
 \node at(3.4,2.5)[ ]{$\delta_2$};
 \node at(3.4,0.95)[ ]{$\delta_3$};
      
        \node at(2.25,0.9)[ ]{$V$};
         \node at(-0.5,-0.2)[right ]{$\Bbb B^2_{\alpha} \times \Bbb B^2_{\beta}$};
          
        \node at(5.5,1.5)[right ]{$\Delta_0=\bigcup_{i=1}^3\delta_i$};
        
  \end{scope}

\end{scope}

 \end{tikzpicture} 
  \caption{$M_0$ and $N_0$ and their normalization}
  \label{fig:trunk}
  \end{figure}
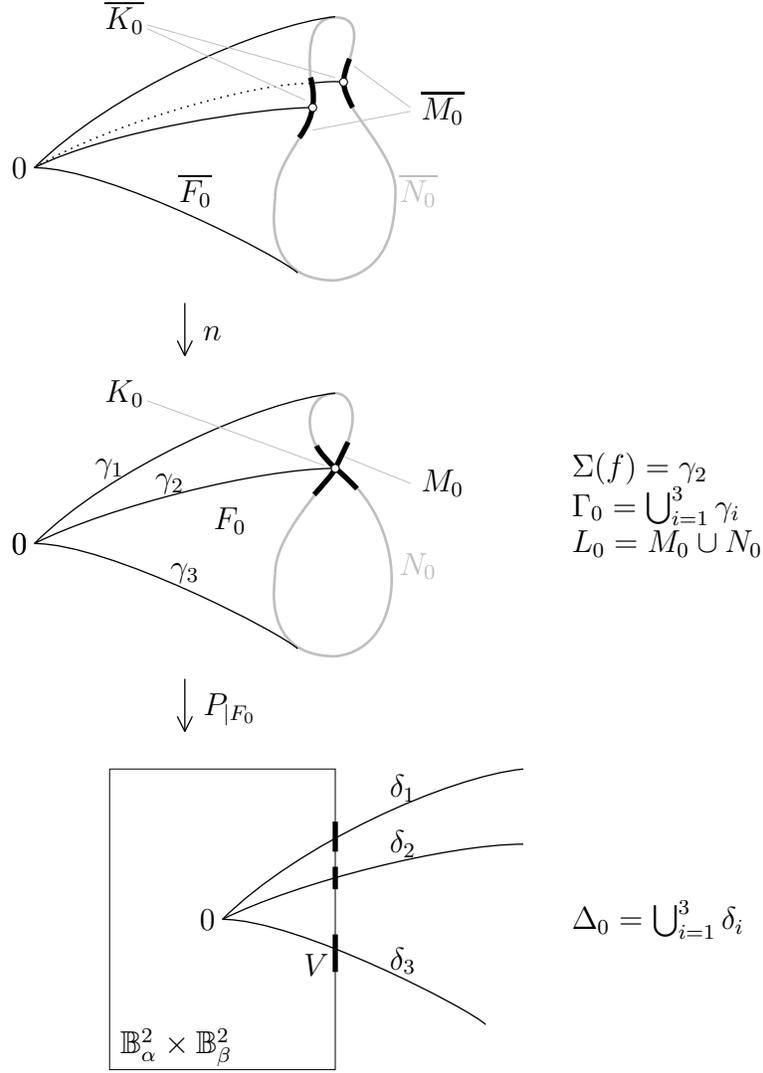

Figure \ref{fig:trunk} shows schematically the construction of $M_0$ and $N_0$ and their normalization when $\Sigma(f)$ is irreducible and $\overline{\Sigma(f)}$ has two  irreducible components.

\begin{proposition}  \label{trunk} The  trunk $N_0$ is a graph manifold with boundary.
  \end{proposition} 
  
\begin{proof} By definition  $\overline{N_0}  =n^{-1} (N_0)$. By construction  $N_0$ does not meet the singular locus $\Sigma (f)$. Therefore the restriction of $n$ on $\overline{N_0}$  is a diffeomorphism from $\overline{N_0}$ to $N_0$. A good resolution  of the pair  $(\overline{F_0},\overline{\Sigma} (f))  $ provides a  decomposition
for $ \overline{L_0}$ as a union of Seifert manifolds such that
$ \overline{K_0}$ is a union of Seifert fibers. As  $ \overline{M_0}$ is 
a tubular neighborhood of
$ \overline{K_0}$ in $ \overline{L_0}$, then the closure $ \overline{N_0}$ of
$( \overline{L_0}  \setminus   \overline{M_0})$ is a   graph manifold with boundary. 
\end{proof}
 
\begin{corollary} \label{sigma} The number of boundary components of the trunk $N_0$  is equal to  the number of irreducible components of the curve  $\overline{\Sigma}(f)$.   
 \end{corollary} 

\begin{proof} In the proof of the above proposition, we show that $N_0$ and $\overline{N_0}$ are diffeomorphic. By construction,  the number of boundary components of the trunk $\overline{N_0}$   equals the number of connected components of $ \overline{K_0} $, which is equal to  the number of irreducible components of the curve  $\overline{\Sigma}(f)$.   
\end{proof}

\begin{rien} \label{Pt}   For each $t \in \Bbb B^2_{\eta} $, the  singular set  $\Gamma_t$  of the restriction of $P$ on $F_t$ is the curve $$ \Gamma _t=  \{ {{\partial f}\over {\partial z}} =0 \}  \cap F_t ,$$ and its discriminant locus is   $\Delta_t=P(\Gamma _t)$.

By  continuity,  we can  choose $\eta$ sufficiently small, $0<\eta <<\alpha$, in such a way that for each $t, |t|\leq \eta$, the  properties  that we already have for $t=0$,  hold for $t \in \Bbb B^2_{\eta} $, {\it i.e.,} 
 
 \begin{itemize}
\item[($1_t$)]    $ L_t \subset \{  |z|  < \gamma  \} $

\item[($2_t$)] $ \Gamma_t $ is a curve   which intersects  ${ \partial B(\alpha)}$ transversally   inside\\ ${ \Bbb S^1_{\alpha} \times int( \Bbb B^2_{\beta}) \times int( \Bbb B^2_{\gamma})}$

\end{itemize}

 \begin{itemize}
\item[($3_t$)]  The restriction of $P$ to $L_t$ is a  finite  cover  with ramification locus   $\Gamma_t \cap ({\Bbb S^1_{\alpha} \times int( \Bbb B^2_{\beta}) \times int( \Bbb B^2_{\gamma})})$ and branching locus $\Delta_t \cap \partial(\Bbb B^2_{\alpha} \times {\Bbb B}^2_{\beta})$.
\item[($4_t$)]  $\Delta_t \cap \partial(\Bbb B^2_{\alpha} \times {\Bbb B}^2_{\beta}) \subset int(V) $   with $V$ as in Proposition \ref{tube}.
\end{itemize}
  \end{rien} 

 \begin{definition} \label{def:vanishingZ} Let $ L( \eta ) =
   f^{-1} (\Bbb B^2_{\eta}  ) \cap   { \partial B(\alpha)}$ and let $M(\eta)$ be the union of the connected components of $L( \eta ) \cap P^{-1}(V)$ which intersect $K_0 $. For any  $t \in \Bbb B^2_{\eta} $,   define: 
   $$M_t := M(\eta) \cap L_t  \hbox{ and } N_t:=\overline{L_t \setminus M_t}.$$ 
   We call $M_t$  and $N_t$ respectively  the {\it  vanishing zone} and  {\it trunk} of $L_t$.
 
 \end{definition}
  
Notice that the choice of $V$  implies that $$M(\eta ) \subset  { \Bbb S^1_{\alpha} \times int( \Bbb B^2_{\beta}) \times int( \Bbb B^2_{\gamma})}.$$
  
 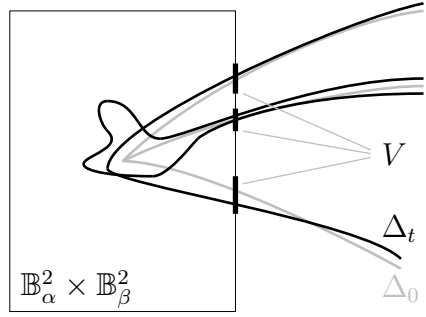
\begin{figure}[ht]
\centering
\begin{tikzpicture} 

% DELTA_0
 \draw[line width=1pt, lightgray  ] 
(1,1.5) .. controls (2,2.5) and (4,3.4)  .. (5,3.5) ;

\draw[line width=1pt, lightgray  ] 
(1,1.5) .. controls (2,2) and (4,2.5)  .. (5,2.5) ;

\draw[line width=1pt, lightgray  ] 
(1,1.5) .. controls (2,1.5) and (4.3,0.3)  .. (4.7,0.07) ;

%DELTA_t

\draw[line width=1pt ] 
(0.8,1.4) .. controls (0.8,2) and (4,3.4)  .. (5,3.6) ;

\draw[line width=1pt ] 
(0.8,1.4) .. controls (0.8,1.2) and (4.2,0.7)  .. (4.7,0.2) ;

%NUMBER 1
\draw[line width=1pt ] 
(1.5,1.8) .. controls (1.8,1.8) and (3.5,2.6)  .. (5,2.6) ;
\draw[line width=1pt ] 
(1.5,1.8) .. controls (1.2,1.8) and (1.2,2.3)  .. (0.8,2.3) ;
 \draw[line width=1pt ] 
(0.8,1.9) .. controls (0.7,2) and (0.6,2.3)  .. (0.8,2.3) ;
 \draw[line width=1pt ] 
(0.8,1.9) .. controls (1,1.7) and (0.6,1.6)  .. (0.5,1.5) ;
 \draw[line width=1pt ] 
(1.4,1.3) .. controls (0.7,1.3) and (0.4,1.4)  .. (0.5,1.5) ;
 \draw[line width=1pt ] 
(2,1.8) .. controls (1.9,1.7) and (1.6,1.3) ..(1.4,1.3) ;
 \draw[line width=1pt ] 
(2,1.8) .. controls (2.1,1.9) and (3,2.4)  .. (5,2.4) ;

% CARRÉ
\draw[line width=0.3pt] 
(-0.5,3.5)--(-0.5,-0.5)--(2.5,-0.5)--(2.5,3.5)--(-0.5,3.5);

% VOISINAGE V
\draw[line width=2pt] 
 (2.5, 2.4)--(2.5,2.8);
 
 \draw[line width=2pt] 
 (2.5, 2.2)--(2.5,1.9);
 
   \draw[line width=2pt] 
 (2.5, 1.3)--(2.5,0.8);

%NODES

\node at(4.3,0.6)[right ]{$\Delta_t $};
\node at(4.3,-0.2)[right, lightgray ]{$\Delta_0 $};

 \node at(-0.5,-0.2)[right ]{$\Bbb B^2_{\alpha} \times \Bbb B^2_{\beta}$};
 %\node at(0.85,1.75)[ lightgray]{$0$};
 
 \node at(4.3,1.6)[right ]{$V $};
  \draw[line width=0.5pt, lightgray] 
 (4.3, 1.65)--(2.6,2.4);

 \draw[line width=0.5pt, lightgray] 
 (4.3, 1.6)--(2.6,1.9);

  \draw[line width=0.5pt, lightgray] 
 (4.3, 1.55)--(2.6,1.2);
 
 \end{tikzpicture} 
  \caption{The curve   $\Delta_t$}
  \label{fig:Deltat}
\end{figure}

\begin{proposition}  \label{eta1} Let $N(\eta )$ be the closure of $L(\eta ) \setminus   M(\eta )$ in $L(\eta )$. There exists a sufficiently small
$\eta$ such that  $f$ restricted to $N( \eta )$ is a fibration on
$\Bbb B_{\eta}^2$.
\end{proposition}
  
 \begin{corollary}\label{eta2} There exists a sufficiently small $\eta $ such that  for all  $ t \in \Bbb B^2_{\eta} \setminus \{ 0 \} $, $N_t$ is orientation preserving  diffeomorphic to $N_0$.   \hfill $\Box$
   \end{corollary}
  
\begin{proof}[Proof of Proposition \ref{eta1}] {We first consider $\eta$ such that  conditions $(1_t)$ to $(4_t)$ of Subsection \ref{Pt} hold for each $t \in \Bbb B^2_{\eta}$.}
  
 i) Let 
   $$\Gamma (\eta)= L (\eta) \cap    \{ {\partial f \over \partial
z} = 0 \}  .$$
 
  Then, the restriction  of $(P,f)$ on  ${ (\Bbb S^1_{\alpha} \times int( \Bbb B^2_{\beta}) \times int( \Bbb B^2_{\gamma}))} \setminus \Gamma (\eta )$ is a submersion.
By $(4_t)$, $\Gamma (\eta )$ does not meet the boundary of $N( \eta )$, hence the restriction of $f$ on the boundary of $N(\eta )$ is a fibration.

\vskip.1in

 ii)  {$(1_t)$ implies that    $$ L(\eta) \subset \partial B(\alpha) \setminus ( \Bbb B_{\alpha}^2 \times
\Bbb B_{\beta}^2 \times \Bbb S_{\gamma}^1)$$
By $(2_t)$, $\Gamma (\eta )$ does not meet $\Bbb B_{\alpha}^2 \times \Bbb S_{\beta}^1 \times  \Bbb B_{\gamma}^2 $, hence the restriction of $f$ on $N(\eta ) \cap (\Bbb B_{\alpha}^2 \times \Bbb S_{\beta}^1 \times  \Bbb B_{\gamma}^2)$ is a fibration. }

\vskip.1in

iii) Now, we have to prove that the restriction of $f $ on $N(\eta )\cap ( \Bbb S_{\alpha}^1 \times
\Bbb B_{\beta}^2 \times \Bbb B_{\gamma}^2)$ is a fibration. Points i) and ii) show that it is a fibration on its  boundary. So, it is sufficient to prove that the projection on the $x$ axis is transverse to $f$ on $N(\eta ) \cap ({ \Bbb S^1_{\alpha} \times int( \Bbb B^2_{\beta}) \times int( \Bbb B^2_{\gamma})})$, {\it i.e.,}   to prove that there exists a sufficiently small $\eta>0$  such that the set
 $$ A  := N(\eta ) \cap ({ \Bbb S^1_{\alpha} \times int( \Bbb B^2_{\beta}) \times int( \Bbb B^2_{\gamma})}) \cap \{ {\partial f \over \partial z} = 0 \}  \cap \{ {\partial f \over \partial y} = 0 \}  $$ is empty. But for a general choice of the coordinates $x$ and $y$,  Subsection \ref{axis1} implies that:
 
 $$ L_0  \cap \{ {\partial f \over \partial z} = 0 \}  \cap \{ {\partial f \over \partial y} = 0 \} = K_0 \subset   int(M_0).$$
 Then, by continuity:
  $$   L(\eta )  \cap \{ {\partial f \over \partial z} = 0 \}  \cap \{ {\partial f \over \partial y} = 0 \}  \subset   int(M(\eta )),$$
    which implies that $A $ is empty. 
 \end{proof}

\begin{rien} \label{vanish} Now,  let us describe more precisely the connected components of  the vanishing zone $M_t$.

 The tubular neighborhood $V$  of $\Delta _0 \cap \partial(\Bbb B^2_{\alpha} \times {\Bbb B}^2_{\beta})$,  used above to obtain the vanishing zone, can be defined   as follows:

 Let $\delta _1,\ldots,\delta _s$ be the irreducible  components of  $\Delta _0$.   Let us fix $i \in \{1,\ldots,s\}$, and  let 

 $$u \mapsto (u^k, \phi _i(u)), \hbox{ where } \phi _i(u) =  \sum_{j=1}^{\infty} a_j u^j$$
be a Puiseux series  expansion of the branch $\delta_i $ of $\Delta_0 $. Let us consider the neighborhood $W_i$ of $\delta_i$ in $\Bbb C^2$ defined by 
$$ W_i =\{ (x,y) \in \Bbb C^2  \ / \ x=u^k, |y - \phi _i(u)|   \leq \theta , u \in \Bbb C \}, $$
where $\theta$ is a positive real number.  

We now choose $\theta$ sufficiently small,  $0<\theta << \alpha$, in such a way that: 

\begin{enumerate}
\item   for each $i=1,\ldots,s$, $W_i$ intersects { $\partial(\Bbb B^2_{\alpha} \times {\Bbb B}^2_{\beta})$} transversally  inside 
$ \Bbb S^1_{\alpha} \times int( {\Bbb B}^2_{\beta}),$
\item the intersection $  V_i = W_i \cap  \partial(\Bbb B^2_{\alpha} \times {\Bbb B}^2_{\beta})$ is a tubular neighborhood of the knot $\delta_i \cap \partial(\Bbb B^2_{\alpha} \times {\Bbb B}^2_{\beta})$,  
\item the  solid tori  $  V_i$ are disjoint.
\end{enumerate}
\vskip0,3cm

{Then $ \bigcup_{i=1}^s  V_i $ can be taken  to be the set $V$ of Proposition \ref{tube}, and we can choose  $\eta << \theta$  such  that conditions $(1_t)$ to $(4_t)$ of Subsection \ref{Pt} hold. In particular, for each $t \in \Bbb B^2_{\eta}$,  we have $  \Delta_t   \cap  \partial(\Bbb B^2_{\alpha} \times {\Bbb B}^2_{\beta})  \subset  int(V)$. }

 Let $\sigma $ be  an  irreducible  component of  $\Sigma(f)$. There exists  $i \in \{1,\ldots,s\}$ such that $P(\sigma )= \delta _i$. We denote by $M (\eta , \sigma )$ the connected component of $P^{-1}(V_i) \cap L(\eta ) $ which  contains the knot   $ K_0 (\sigma) = \sigma \cap  {\partial B(\alpha)}$ of $\sigma $ in $ {\partial B(\alpha)}$. 
 
 By definition,  the $3$-dimensional manifold $M_t(\sigma) =  M( \eta , \sigma ) \cap L_t$ is connected, and we  obtain: $$ M_t = \bigcup_{j=1}^r M_t(\sigma_j),$$
 where $\{ \sigma_j , 1\leq j \leq r \} $ is the set of the irreducible components of $\Sigma (f)$.
 
For each $j=1,\ldots,r$, let  $\bar{r}_j$ be the number of irreducible components of the curve $n^{-1}(\sigma_j)$. The boundary of $M_t(\sigma_j)$ consists of $\bar{r}_j$ tori.
\end{rien} 

\begin{definition} \label{def:vanishZ}
$M_t(\sigma)$   is the {\it vanishing zone of $L_t$ along} $\sigma$.
\end{definition}

  Proposition \ref{trunk}, Corollary \ref{eta2} and the construction in Subsection \ref{vanish} summarize in  the following theorem:

 \begin{theorem} \label{main1} 
 
  \begin{enumerate}  
 \item The boundary $L_t$ of the Milnor fiber of $f $ decomposes as the union 
 $$L_t = N_t \cup M_t,$$

\item $N_t \cap M_t$  is a disjoint union of $r$ tori, where $r$ is the number of irreducible components of the curve $\overline{\Sigma}(f)$,  
\item $N_t$ is a     graph manifold orientation preserving  diffeomorphic   to the trunk $\overline{N_0}$,
\item Let $\sigma_1,\ldots,\sigma_r$ be the  irreducible components of $\Sigma(f)$.   The  connected components of the vanishing zone $M_t $ are the manifolds $M_t(\sigma_j), j = 1 \ldots r$. \hfill $\Box$
\end{enumerate}
\end{theorem}

\begin{corollary}\label{connected} The manifold $L_t$ is connected. 
\end{corollary}

\begin{proof}  The number of connected components of
$\overline {F_0}$ and $\overline {L_0}$ is equal to the number of irreducible
components of $f$. The intersection between two irreducible components of
$f=0$ furnishes at  least one irreducible component of the singular locus
$\Sigma(f)$ and a corresponding connected component of the vanishing
zone. Hence, the constructions given here  show that after the gluing of
all connected components of the vanishing zone with the trunk, we obtain
a connected manifold $L_t$. 
\end{proof}

\begin{remark} Corollary \ref{connected} implies that the Milnor fiber 
$F_t$ is connected.  In fact, as the singular locus of $f$ has dimension $1$, $F_t$
is connected by a much more general result of  M. Kato and  Y. Matsumoto
in \cite{K-M}.
\end{remark}

\begin{remark} \label{remark1}To prove that $L_t$ is a  graph manifold, we still have to prove that  $M_t (\sigma )$ is a   graph manifold for any irreducible component $\sigma $ of $\Sigma (f)$. This will be done in Section 4.
\end{remark}

{\begin{rien} \label{reduction} {\bf Reduction  to a smooth branch of $\mathbf \Gamma_t$}

 Let us fix a branch $\sigma$ of $\Sigma(f)$ and let  $$ u  \mapsto (u^k,\phi(u), \psi(u) )$$ 
be a Puiseux parametrization of $\sigma$. 

 Let us consider the analytic morphism $\Theta :\Bbb C^3 \rightarrow \Bbb C^3$ defined by 
$$\Theta(x,y,z)=(x^k,y+\phi(x), z + \psi(x)).$$
   \vskip.1in
    Let $g : (\Bbb C^3,0) \longrightarrow (\Bbb C,0)$ be the composition $g = f \circ \Theta$. Then  $ \sigma' = \Theta^{-1}(\sigma)$ is the $x$-axis. Moreover, a direct computation of the derivative of $g$ shows that $\sigma'$ is a branch of the singular locus of $g$.
    
      Let $M_t(f,\sigma)$ (resp.  $M_t(g,\sigma')$)  be the vanishing zone of $f$  along $\sigma$ (resp. of $g$ along $\sigma'$) defined in the boundary of the polydisc $B(\alpha)$ (resp. $B(\alpha^{1/k})$)  as in \ref{vanish}. The construction given in \ref{vanish} leads directly to: 
      
\begin{lemma} {The  restriction $$\Theta_{\mid M_t(g,\sigma')} : M_t(g,\sigma') \rightarrow M_t(f,\sigma)$$ is a diffeomorphism. \hfill $\Box$}
\end{lemma} 
   
\noindent
 In the sequel, we assume that $\sigma$ is the $x$-axis. In particular, the vanishing zone $M_t(\sigma)$ along $\sigma$ is nothing but 
$$M_t(\sigma) =L_t \cap  ( \Bbb S^1_{\alpha} \times  \Bbb B^2_{\theta} \times  \Bbb B^2_{\gamma}), \ 0<\eta << \theta << \alpha .$$
 \end{rien}
}

\section{A parametrization theorem}\label{sec:parametrization}

  This section can be read independently of the others. It contains  the {\it Parametrization Theorem} \ref{par} which will play a key role  in Section \ref{proof} to perform the carrousel in family. To prove it, we construct an embedded  resolution for a family of plane curve germs parametrized  by a punctured disc. This provides a {\it Resolution in family Theorem} \ref{th:resolution} that we state and prove at the end of this section. This latter independent result will not be used in the rest of the paper.
  
We consider a reduced   holomorphic  germ $h : (\Bbb C^3,0) \rightarrow (\Bbb C,0) $ such that $ h(x,0,0)=0$ for all $ x \in \Bbb C$. Let  $(H,0)$ be the  hypersurface  germ with equation $h=0$. 

 Let  $\alpha >0 $  be such that the polydisc $B(\alpha )$   is a
Milnor polydisc for $h$ (Definition \ref {def:polydisc}). For each $x \in  \Bbb B_{\alpha }^2 \setminus   \{0\}$, we denote by $h_x : (\Bbb C^2,0) \rightarrow (\Bbb C,0)$ the germ defined by: 
$h_x(y,z)= h(x,y,z)$. Hence  the germ $h_x$  is either  non-singular  or has an isolated singular point at $(x,0,0)$  for all $x \in \Bbb B_{\alpha  }^2 \setminus   \{0\} $. For each $\alpha' , \ 0<\alpha'  \leq \alpha $, there exists $\epsilon>0$, ($\epsilon  << \alpha' $), such that for each $x \in  \Bbb S^1_{\alpha' } $, $\{x\} \times \Bbb B^4_{\epsilon }$ is a Milnor ball for  $h_x$.

\begin{definition} \label{sh} {Let $\alpha' $ such that  $0 < \alpha'  \leq \alpha $} and let $\epsilon$  be associated with $\alpha ' $ as above. A {\it {sheet}  of $H$ along 
$ \Bbb S^1_{\alpha' } \times  \{ 0 \} \times  \{ 0 \} $ } is  the closure of a connected component of the intersection $H \cap (  \Bbb S^1_{\alpha' }  \times ( \Bbb B^4_{\epsilon} \setminus \{0\})) $. 
\end{definition}

\begin{theorem}[{\it Parametrization Theorem}]    \label{par}   {There exists  $\alpha >0$  sufficiently small such that for each $\alpha' $ with $0<\alpha'  \leq \alpha $,  the following holds.} Let  $G$ be a sheet         of $H$ along $\Bbb S^1_{\alpha' }  $.   There  exist $d, i$ and $j \in \Bbb N^*$, and two convergent power series $b(x^{1/d},u) \in \Bbb C\{x^{1/d}\}\{u\}$ and $c(x^{1/d},u) \in \Bbb C\{x^{1/d}\}\{u\}$    with $b(x^{1/d},0) \neq 0$ and $c(x^{1/d},0)\neq 0 $, such that 

 $$ (s,u)   \longmapsto (s^d,u^i b(s,u), u^{j} c(s,u)) $$
  with $s \in \Bbb S^1_{{\alpha'}^{1/d}}$ is a parametrization of $G$.  

\end{theorem}  

\begin{proof}[Beginning of proof of Theorem  \ref{par}] 
 
  We start by taking $\alpha $ such that $B(\alpha )$ is a Milnor polydisc for $h$. We will later decrease $\alpha$ (if necessary) in order that three extra
conditions hold: Condition \eqref{C1} stated below and Conditions \eqref{C2} and \eqref{C3}  introduced in the proof of Lemma \ref{LemPar}.

 Let us write $h(x,y,z)$ as the sum 
$$ h(x,y,z) = \sum_{n=0}^{\infty} c_n (x,y,z),$$
where for all $  n \in \Bbb N$,  
$$c_n(x,y,z) = \sum_{k=0}^n c_{n,k}(x) y^kz^{n-k},$$
with $c_{n,k}(x) \in \Bbb C\{x\}$. 

Let $m$ be the least integer such that $c_m(x,y,z) \neq 0$.  Perhaps after performing a change of variables, one can assume that $c_{m,0}(x) \neq 0$   in $\Bbb C\{x\}$.

{By  decreasing   $\alpha $ if necessary,  we can assume that the following condition \eqref{C1} holds:

\[
 \hbox{For all  $s \in  \Bbb B^2_{\alpha } \setminus   \{ 0 \}, \ c_{m,0}(s) \neq 0$}  \tag{C1}\label{C1}
\]
}

   Let us  treat first the case $m=1$, {\it i.e.,} $h_x$ is non-singular.  If $m=1$, there exists $c(x,z)\in  \Bbb C\{x,z\}$ such that  $h(x,0,z)=c(x,z)z$.   Since $c(x,0)=c_{m,0}(x)$, Condition \eqref{C1} implies that   in $B(\alpha )\setminus \{x=0\}$,  we can change the coordinate $y$ in $y=c(x,z)y'$. There exists $g(x,y',z)\in  \Bbb C\{x,y',z\}$ such that:
   $$h(x,y,z)=c(x,z)(y'g(x,y',z)+z).$$

Let us consider     $F(y',z) =   y'  g( x , y',z) +   z $ as an element of $A\{y',z\}$ where $A = \Bbb C\{ x \}$. As $F(0,z)=z$, we can apply  the Weierstrass preparation theorem  (for example see \cite {Z-S}, vol.2, p.139-141)  to obtain  $R( x,y') \in \Bbb C\{ x \}\{y'\}$ such that :
$$F(y',z)=0 \Leftrightarrow z = R(x ,y').$$

  {Then for each $\alpha'   $ such that $0<\alpha'  \leq \alpha $,} the intersection  $H \cap (  \Bbb S^1_{\alpha' }  \times ( \Bbb B^4_{\epsilon} \setminus \{0\})) $ is  connected and parametrized by $(x,y',R(x,y'))$.   This achieves the proof when $m=1$.

  When $m>1$, we need Lemma \ref{LemPar}   to prove  Theorem \ref{par}.  To state Lemma \ref{LemPar}, we consider, for each $x \in  \Bbb B^2_{\alpha } \setminus \{0\}$, the minimal good resolution   $\pi_x : Y_x \rightarrow \{x\} \times  \Bbb B^4_{\epsilon }$   of $h_x$, {\it i.e.,} the minimal   composition of blow-ups of points such that the curve $(h_x \circ \pi_x)^{-1}(0)$ is a normal crossing divisor.  We denote by $E_x= \pi_x^{-1}(x,0,0)$ the exceptional divisor of $\pi_x$.

\begin{lemma}  \label{LemPar}   We assume $m>1$. For a chosen  $x \in  \Bbb B_{\alpha }^2 \setminus   \{0\}$,  let  $h_{1,x}$ be an irreducible  factor of $h_x$,  let $\tilde h_{1,x}$ be its strict transform by  $\pi_x$ and let $P=E_x \cap \tilde h_{1,x}$.  For  a sufficiently small $\alpha $,  we can choose  coordinates $(u,v)$ at $P$ in $Y_x$ such that: 
\begin{enumerate}

\item\label{eq:lempar1}  $u=0$ is a local equation for $E_x$ in $Y_x$. 

\item\label{eq:lempar2} There exist  
three   integers $d,i,j $ in $ \Bbb N^*$,  two polynomials 
$\phi(x^{1/d},u,v)$ and $\psi(x^{1/d},u,v) $ in  $ \Bbb C\{x^{1/d}\} [u,v]$, where  $\phi( x^{1/d},0,v)$ and $\psi( x^{1/d} ,0,v) $ are not identically $0$,  and   $s  \in   \Bbb B^2_{{\alpha }^{1/d}} \setminus \{ 0\} $ with $s^d=x$  such that:   

$$ \pi _x ( u,v) = (s^d,u^i\phi(s,u,v), u^j \psi (s,u,v)).$$
\item\label{eq:lempar3} There exist an integer $M \in \Bbb N^*$ and two convergent power series   $c(x^{1/d}) \in ( \Bbb C\{x^{1/d}\} \setminus \{ 0 \} )$ and $g(x^{1/d},u,v) \in \Bbb C \{x^{1/d}\}\{u,v\}$ such that, for the value $s$ defined just above, we have:
$$(h \circ \pi_x)(u,v) = u^M \big(  u g(s,u,v) + c(s) v \big).$$
\end{enumerate}
\end{lemma}

\begin{proof}[Proof of Lemma \ref{LemPar}]
 
We  start with the blow-up $\pi_{1,x}$ of $(x,0,0)$ in $\Bbb C^2$, {\it i.e.,}
$$ \pi_{1,x}  : Y_{1,x} \rightarrow \{ x \} \times \Bbb B^4_{\epsilon }. $$

Let $E_{1,x}= (\pi_{1,x})^{-1}(x,0,0)$ be the exceptional divisor of $\pi_{1,x}$. As $c_{m,0}(x) \neq 0$, the  line $y=0$ is not  tangent to the curve $h_x=0$. We  will write  the intersection points  between   $ E_{1,x}$ and the strict transform of $ h_x $,   with the help of coordinates  $(u_1,v_1)$  given by  the standard chart on  $(\pi_{1,x})^{-1 }  ( \{ x \}   \times \Bbb B^4_{\epsilon }  )$ defined by 
$$  \pi_{1,x}(u_1,v_1) = (x,u_1,u_1v_1).$$

In the local coordinates $(u_1,v_1)$, we have: 

$$(h_x \circ \pi_{1,x} )(u_1,v_1) = u_1^m \bigg( \sum_{k=0}^m c_{m,k}(x) v_1^{m-k }+  u_1 g_1 (x,u_1,v_1)  \bigg) \hskip1cm$$
where 
$$ g_1(x,u_1,v_1) = \sum_{m'=m+1}^\infty u_1^{m'-m-1} c_ {m'}(x,1,v_1).$$
Then  the intersection points  between   $ E_{1,x}$ and the strict transform of $ h_x $ are the points $(x,0,v_1)$ such that $v_1$ is a root of the polynomial
$$ Q(v_1) =  \sum_{k=0}^m c_{m,k}(x) v_1^{m-k } \in {\Bbb C}\{x\} [v_1].$$

There exists  an integer $e>0$ such that the decomposition  field of the polynomial $Q$ is the fraction field  $\Bbb C \{\{  x^{1/{e}}\}\}$ of $\Bbb C \{  x^{1/{e}}\}$ (for example see D.Eisenbud \cite {Ei}, p.295).  Let $P_1$ be the intersection point between $E_{1,x}$ and  the strict transform of $h_{1,x} $ (by $\pi_{1,x}$).    There exists  a   root $r_1 \in  \Bbb C \{\{  x^{1/{d_1}}\}\}$  of $Q$, where   $d_1\leq e$ is the  minimal integer such that $r_1  \in \Bbb C \{\{  x^{1/{d_1}}\}\}$, and a complex number $s_1$ which satisfies  $s_1^{d_1}=x$,  such that  $ P_1=(0,r_1(s_1))$.     Let $\delta $ be a $d_1-$th root of the unity  ($\delta ^{d_1}=1$). The strict transform  of $h_x$  meets  also $E_{1,x}$ at the $d_1$ distinct points  $(0,r_1( \delta s_1))$.

{The Galois group of  $\Bbb C \{\{  x^{1/{e}}\}\}$ acts on the roots of $Q$. The orbit of $r_1$ is the set $\lbrace  r_{\delta} , \ \delta^{d_1}=1\rbrace $,  where  $ r_{\delta}( x^{1/{d_1}}) = r_1( \delta  x^{1/{d_1}} ).$  As  $h_{1,x} $   is irreducible, the root $r_1$ is unique up to  the action of the  Galois group.

We proceed as above with chosen  representatives of the other  orbits of the roots of $Q$. Thus,  we find     all the intersection points of the strict transform of $h_x$ (by $\pi _{1,x}$) with $E_{1,x} $. The  map $\pi_{2,x}$ is the  blow-ups of all  these intersection points.

{By  decreasing again   $\alpha $ if necessary,  we can assume that the following condition \eqref{C2} holds for each pair $r$ and $r'$ of  roots of $Q$ which represent   distinct orbits under the action of the Galois group:

\[
 \hbox{ $r \big( \Bbb B^2_{{\alpha }^{1/e} } \setminus   \{ 0 \}\big) \cap r' \big( \Bbb B^2_{{\alpha }^{1/e} } \setminus   \{ 0 \}\big) = \emptyset $}  \tag{C2}\label{C2}
\]
 }}

As $\Bbb C \{\{  x^{1/{d_1}}\}\}$ is nothing but the field of convergent Laurent power series in the variable $x^{1/{d_1}}$, there exists $l_1 \in \Bbb N^*$ such that 
$$x^{l_1} r_1(x^{1/{d_1}}) \in \Bbb C\{x^{1/d_1}\}$$ 

We  consider new local coordinates $(\tilde u_1, \tilde v_1)$ in $Y_{1,x}$ centered at $( 0,r_1(x^{1/d_1}))$ by setting: 
$$  ( u_1,  v_1 )= (x^{l_1} \tilde u_1, \tilde v_1 + r_1(x^{1/d_1})) $$

We then have:  
$$ \pi_{1,x} (\tilde u_1,\tilde v_1) = (x,u_1,u_1v_1)=(x,x^{l_1} \tilde u_1, x^{l_1} \tilde u_1( \tilde v_1 + r_1(x^{1/d_1}))) \hskip1cm (*)$$

As $x^{l_1} $ and $ x^{l_1} ( \tilde v_1 + r_1(x^{1/d_1}))$ are   in 
 $ \Bbb C\{x^{1/d_1}\} [\tilde u_1, \tilde v_1] $ and as $\tilde u_1=0 $ is the local equation of $E_{1,x}$ at the point  $P_1$, statements (1) and (2) of Lemma \ref{LemPar} are proved for $\pi _{1,x}$.

  \vskip.1in
  When we  perform   $\pi_{2,x}$, we  blow up  $P_1$ in $Y_{1,x}$. In order to write $\pi _{2,x}$ in one of the  two standard charts around $(\pi _{2,x})^{-1}(P_1)$, we perform in $(\ast)$ one of the two following substitutions: $ ( \tilde u_1,   \tilde v_1) = (u_2, u_2v_2)$ or $ ( \tilde u_1,   \tilde v_1) = (u_2 v_2, v_2)$. If necessary, we follow it by a  new change of coordinates of the type:
   $$  ( u_2,  v_2 )= (x^{l_2} \tilde u_2, \tilde v_2 + r_2(x^{1/d_2})), $$
   where $x^{l_2} r_2(x^{1/{d_2}}) \in \Bbb C\{x^{1/d_2}\}$ is    defined as before. 

    Then, points \eqref{eq:lempar1} and \eqref{eq:lempar2} of Lemma \ref{LemPar} are also proved for $\pi _{2,x}  \circ   \pi _{1,x}$. By  finite iteration  (where we take care of  conditions \eqref{C1} and  \eqref{C2}  to shrink  $\alpha $ if necessary),  they are also proved for $\pi _x$. As $\pi _x$ is a good resolution of $h_x$, the strict transform $\tilde h_x$ is transverse to $E_x$ at $P$ and has multiplicity $1$.  By a  direct computation of $ h_x  \circ  \pi_x$, with the help of point \eqref{eq:lempar2} of Lemma \ref{LemPar}, we obtain point \eqref{eq:lempar3} of Lemma \ref{LemPar}:
    $$(h \circ \pi_x)(u,v) = u^M \big(  u g(s,u,v) + c(s) v \big),$$
    where $c(s)\in   \Bbb C\{s^{1/d}\} \setminus \{ 0 \}. $

{By  decreasing    $\alpha $  again if necessary,  we can assume that the following condition \eqref{C3} holds:

 \[
 \hbox{For all $s \in  \Bbb B^2_{{\alpha }^{1/d} } \setminus   \{ 0 \},\ c(s) \neq 0$}  \tag{C3}\label{C3}
\]
Then, for each $x \in \Bbb B_{\alpha }^2 \setminus \{0\}$, $\pi_x$ is the minimal resolution of $h_x$.}
 \end{proof}

  \begin{rien}{\bf Resolution in family}\label{resfa}

   {Let $\cal U(\alpha )$ be the set $ \{ x \in \Bbb C , \ 0<\vert x \vert  <\alpha \}$. If $B(\alpha )$ is a Milnor polydisc for $h$ and if $\alpha  $ satisfies conditions \eqref{C1}, \eqref{C2} and \eqref{C3},  then  the proof of Lemma \ref{LemPar}   shows that   the dual resolution graph of  $\pi_x $  does not depend on $x$ and then that  the family of  germs $h_x$,  parametrized by $x\in \cal U(\alpha )$, has  constant topological type. To state a {\it resolution in family} theorem  in our particular case (our family is parametrized by $x$ where  $x\in \cal U(\alpha )$), we use the following definition:  

\begin{definition} A  {\it punctured disc} is a $1$-dimensional analytic open manifold isomorphic to the standard punctured disc  $\cal U=\lbrace  s\in \Bbb C, \ 0< |s|  <1 \rbrace$.
\end{definition}

  In the sequel, we will consider blow-ups of  analytic manifolds along punctured discs in the sense of \cite[page 378]{Lo}. 

\begin{lemma}    \label{res}   Assume  $\alpha >0$ is  such that  $B(\alpha )$ is a Milnor polydisc for $h$ and conditions \eqref{C1}, \eqref{C2} and \eqref{C3} are satisfied. Then there exist a smooth $3$-dimensional open analytic manifold  $Y$ and a composition of a finite number of  blow-ups  along  punctured discs $\pi : Y \ \rightarrow int(B(\alpha ) \setminus   \{x=0\})$  such that:

\begin{enumerate}
\item \label{res1} $\pi$ restricted on $Y\setminus \pi ^{-1} (\cal U(\alpha )\times \{0\}\times \{0\})$ is an isomorphism.
\item  \label{res2} $E=\pi ^{-1} (\cal U(\alpha )\times \{0\}\times \{0\})$ is an analytic normal crossing divisor.
\item \label{res3} For each $x\in\cal U(\alpha )$, there exists $\epsilon>0$  such that the closure of $ \pi ^{-1} (H\cap \{x\}\times (\Bbb B^4_{\epsilon}  \setminus \{0\}))$ in $Y$  is a disjoint union of smooth discs which cut  transversally  the exceptional divisor $E$  at smooth points.
\end{enumerate} 
\end{lemma} 
}
\begin{proof}
Let  $ \pi_{1} : Y_1 \rightarrow   int(B(\alpha  )\setminus   \{ x=0 \}) $  be  the blow-up  {of  $ int(B(\alpha  )\setminus   \{ x=0 \})$  along $\cal U(\alpha )\times   \{ 0 \} \times   \{ 0 \}$}. Let $E_1$ be the exceptional divisor of $\pi_1$ and $H_1$ be the strict transform of   $H\cap int(B(\alpha  )\setminus   \{ x=0 \}) $. We have:
   $$ \pi_1^{-1}(H\cap int(B(\alpha )\setminus   \{x= 0 \}) )=E_1\cup H_1.$$

  We follow the steps of the proof of Lemma \ref{LemPar}, but here we do it for all $x\in \cal U(\alpha )$. Condition \eqref{C1} implies that we can write $\pi_1$ in the same chart  for all 
  $x\in  \cal U(\alpha )$.  More precisely, in the chart $(u_1, v_1)$ used in the proof of  Lemma \ref{LemPar}, the  connected  components  of  $ H_1\cap E_1$ are  parametrized by $\{ (s ,0,r(s)), \ s^{e} \in \cal U(\alpha)\},  $ where $r$ are roots of $Q$ {which represent  all the distinct orbits   under the action of the Galois group of $\Bbb C\{\{x^{1/e}\}\}$. } Condition \eqref{C2} implies that $H_1\cap E_1$ is a disjoint union of punctured discs. Let $\pi _2$ be the blow-up  of  {$Y_1$ along $E_1 \cap  H_1$}. We iterate the  process to obtain $\pi = \pi_k \circ...\circ  \pi_2 \circ \pi_1 $.
 {For each  $ x\in \cal U(\alpha ) $, the restriction of $\pi$  on $\pi ^{-1}( \{x\}\times \Bbb B^4_{\epsilon} )$}   is the minimal resolution $\pi_{x}$ of $h_x$. The construction of $\pi$ implies \eqref{res1} and   \eqref{res2}. Condition \eqref{C3} implies  \eqref{res3}.   
\end{proof}

 \end{rien} 

\noindent
{\it End of proof of Theorem \ref{par}.}   Let $\alpha >0$ be as in Lemma \ref{res}. Let us fix $\alpha' >0$ such  that $0<\alpha' < \alpha $ and let  $\epsilon>0$ ($\epsilon  << \alpha ' $) be such that  $\{x\} \times \Bbb B^4_{\epsilon }$ is a Milnor ball for  $h_x$  for each $x\in   \Bbb S^1_{\alpha' } $.  Let $\tilde G$ be a connected component  of the strict transform  of $H\cap (\Bbb S^1_{\alpha' } \times \Bbb B^4_{\epsilon })$ by $\pi$.   Then  $G= \pi (\tilde G)$ is a sheet of $H$ along $ \Bbb S^1_{\alpha' }\times  \{ 0 \} \times  \{ 0 \} $ (Definition \ref{sh}). {Lemma \ref{res}  implies that $\tilde G$ is a solid torus such that  for each  $x\in \Bbb S^1_{\alpha' }$, the strict transform  of $h_x$ is a disjoint union of meridian discs of $\tilde G$.}

  We  choose   the  sheet  $G$ of $H$   which contains the irreducible factor  $h_{1,x}$ of $h_x$ considered in   Lemma \ref{LemPar}. We now consider $x$ as a variable in $ \Bbb S^1_{\alpha' }$. Point \eqref{eq:lempar3}  of  Lemma \ref{LemPar}   implies that, in a neighborhood  of  $\tilde G \cap E$ in $Y$,  there are  local coordinates $(s, u,v)$ such that   $s^{d}=x$ and there exist  a positive integer $M$,  $c(s)  \in ( \Bbb C\{ s \} \setminus \{ 0 \})$ and $g(s ,u,v) \in \Bbb C \{ s, u,v\}$  which satisfy:  
$$(h \circ \pi )(s,u,v) = u^M \big(  u g(s,u,v) + c(s) v \big).$$

To parametrize $G$, we have to solve the  equation:
$$   u g(s,u,v) + c(s) v =0.$$
But if $s^{d} \in  \Bbb S^1_{\alpha' }$,  we have $c(s)\neq 0$ by Condition \eqref{C3}. Let us perform  the change of coordinate $u'=   c(s)^{-1} u$. We  obtain:
$$(h \circ \pi )(s, u',v) =  {u'}^M c(s)^{M+1} \big(  u' \  g(s ,  u'c(s) ,v) +   v \big).$$
We replace $u'$ by $u$. Now the equation of $\tilde G$ is given by:
$$  u \  g(s ,  { uc(s)},v) +   v  = 0.$$

Let us consider   $F(u,v) =   u \   g( s ,  { uc(s)},v) +   v  = 0$ as an element of $A\{u,v\}$ where $A = \Bbb C\{ s \}$. As $F(0,v)=v$, we can apply  the Weierstrass preparation theorem (for example see \cite {Z-S}, vol.2, p.139-141)  to obtain  $R( s ,u) \in \Bbb C\{ s \}\{u\}$ such that 
$$F(u,v)=0 \Leftrightarrow v = R(s ,u).$$

This leads to:
$$h \circ \pi   ( s,u,R(s,u) )=0. $$

This equality and  point \eqref{eq:lempar2} of Lemma \ref{LemPar}  imply that $h$ vanishes on  $ \{  (s^d,u^i\phi(s,u,R (s,u)), u^j \psi (s,u, R(s,u))) ,  u \in \Bbb B^2_{ \epsilon }  \} $.

For each $s^{d} \in  S^1_{\alpha' }$,  we set  $ b(s,u) =   \phi(s,u,R(s,u))$ and $  c(s,u)=  \psi(s,u,R(s,u))$. We have a parametrization $$    \Bbb S^1_{\alpha'^{1/d}}  \times \Bbb B^2_{\epsilon } \rightarrow  G $$ given by 
$$ (s,u)   \longmapsto (s^d,u^i b(s,u), u^{j} c(s,u)).$$
 \end{proof}
 
 \begin{rien}{\bf Resolution in family: general case}

{Let us  consider a reduced  holomorphic  germ $f : (\Bbb C^3,0) \rightarrow (\Bbb C,0)$. We choose generic coordinates as described in Subsection \ref{axis1}.  Let  $\alpha >0 $  be such that the polydisc $B(\alpha )$   is a
Milnor polydisc for $f$ (Definition \ref{def:polydisc}). Let $\sigma$ be an irreducible component of the singular locus of $f$ and let $ \sigma ^{*}=\sigma \cap ( int(B(\alpha ) \setminus   \{x=0\})).$ 

\begin{theorem} \label{th:resolution}

 There exists a sufficiently small  $\alpha >0$, an  open analytic manifold $ V,$ neighborhood of  $ \sigma ^{*}$ in $ int(B(\alpha) \setminus   \{x=0\}),$  and a composition of a finite number of  blow-ups  along  punctured discs $\pi :  \tilde V \ \rightarrow   V$ {starting with the blow-up of $V$ along $\sigma^*$}  such that:
\begin{enumerate}
\item \label{res1} $\pi$ restricted on $\tilde{V} \setminus \pi ^{-1} (\sigma ^{*})$ is an isomorphism.
\item  \label{res2} $E=\pi ^{-1} (\sigma ^{*})$ is an analytic normal crossing divisor.
\item \label{res3} For each $p=(x,y,z)\in \sigma ^*$,   $\pi$ restricted on    $ \pi ^{-1}(V \cap ( \{ x \} \times \Bbb C^2) ) $  is an embedded resolution of  the plane  curve germ, at $p$,  $ V \cap ( \{ x \} \times \Bbb C^2) \cap \{ f=0\}$.

\end{enumerate}

  \end{theorem}

\begin{proof}
  By   Subsection  \ref{reduction}, we can assume that the singular curve  $\sigma $ is the $x-$axis. Then,  we take $\alpha $ as in Lemma  \ref{res}.  For each $\alpha' , \ 0<\alpha'  < \alpha$, there exists $\epsilon (\alpha'  )>0$, {which depends analytically on $\alpha ', $} such that,  for each $x \in  \Bbb S^1_{\alpha' } $,  $\{x\} \times \Bbb B^4_{\epsilon (\alpha'  ) }$ is a Milnor ball for  $ \{ f=0\} \cap  ( \{ x \} \times \Bbb C^2) $.
    {Then we   take}   
   $$V= \bigcup _{0<\alpha' < \alpha } (  \bigcup _{x\in  \Bbb S^1_{\alpha ' } } \ \ int ( \{ x\} \times  \Bbb B^4_{\epsilon (\alpha ' )})),$$  and we apply  Lemma \ref{res}.

\end{proof}
}
\end{rien}

\section{${\mathbf M_{\mathbf t}}$ is   a graph manifold:   the proof} \label{proof}

{The aim of this   section  is to prove the main result of this paper:  

\begin{theorem}\label{main}  There exists a sufficiently small $\eta $ such that for all $t$ with  $0<  |t|  \leq \eta $,  $M_t$ is  a       graph  manifold whose Seifert pieces have oriented basis.
\end{theorem} 

According to Theorem \ref{main1},  we have to prove that for each irreducible component   $\sigma$  of the singular locus $\Sigma(f)$, the vanishing zone  $M_t(\sigma)$  of $L_t$ along $\sigma $ (Definition \ref{def:vanishZ}) is a  graph manifold. 
 Using Subsection \ref{reduction}, we can assume that $\sigma$ is the $x$-axis. Then,  $M_t(\sigma)$  is given by:
$$M_t(\sigma) =L_t \cap  ( \Bbb S^1_{\alpha} \times  \Bbb B^2_{\theta} \times  \Bbb B^2_{\eta}), \ 0<\eta << \theta << \alpha .$$

 As in Subsection \ref{Pt}, we consider  the projection $P\colon \Bbb C^3 \to \Bbb C^2$ defined by $P(x,y,z)=(x,y)$ and its restriction to the Milnor fiber $F_t= f^{-1}(t) \cap (\Bbb B^2_{\alpha} \times \Bbb B^2_{\beta} \times \Bbb B^2_{\gamma})$.  The  critical locus of $P_{\mid F_t}$ is  $\Gamma_t =  \{\frac{\partial f}{\partial z}=0\} \cap F_t$ and  the set of its critical values is  $\Delta_t= P(\Gamma_t)$. 
   
 As $\sigma$ is the $x$-axis, the image of   $M_t(\sigma)$ by $P$  is the solid torus  $ \Bbb S^1_{\alpha } \times \Bbb B^2_{\theta}$. The restriction of $P$ to $M_t(\sigma)$ is a finite  cover  which is ramified over the braid $\Delta_t \cap   (\Bbb S^1_{\alpha } \times \Bbb B^2_{\theta}).$ But,  if $t\neq 0$,   $\Delta_t$ is not a representative of a curve germ.  In particular,  we cannot use  Puiseux series  expansions to describe  $\Delta_t \cap   (\Bbb S^1_{\alpha } \times \Bbb B^2_{\theta})$ 
 as an iterated torus link.  In order to circumvent  this difficulty,  we consider the germ  
   $$\Psi : (\Bbb C^3,0) \rightarrow (\Bbb C^3,0)$$  defined by $$\Psi(x,y,z) = (x,y,f(x,y,z)).$$  
   The critical locus of $\Psi  $ is the germ of surface 
   $$H := \{ {{\partial f}\over {\partial z}} =0  \}, $$  
   and its discriminant locus is the image $$H' = \Psi (H).$$

        We will study the restriction of $\Psi$ on the manifold 
   $M(\eta,\sigma)$  which is the  union of the $M_t (\sigma )$ where $t \in \Bbb B^2_{\eta }$ (see Subsection  \ref{vanish}).  It is a finite  cover over the product $\Bbb S^1_{\alpha } \times \Bbb B^2_{\theta }\times \Bbb B^2_{\eta }$ which is ramified over the intersection  $ H'\cap (\Bbb S^1_{\alpha } \times \Bbb B^2_{\theta }\times \Bbb B^2_{\eta }) $.

 On one hand, let us fix $t \in \Bbb S^1_{\eta}$ and let  $\Psi_t$ be the restriction of $\Psi$ on  $ M_t (\sigma)$:
   $$\Psi_t \colon M_t(\sigma) \rightarrow \Bbb S^1_{\alpha} \times \Bbb B^2_{\theta}\times \{t\}.$$ 
  By definition $\Delta_t \times \{t\}= H' \cap (\Bbb B^2_{\alpha} \times \Bbb B^2_{\theta}\times \{t\}).$
 So, the set of the ramification values  of the finite cover $\Psi_t$ is the braid 
  $$H'_t = H'\cap (S^1_{\alpha} \times \Bbb B^2_{\theta}\times\{t\}).$$

      On the other hand, let us fix $a \in \Bbb S^1_{\alpha}$ and let us consider  $M^{(a)}(\sigma)= M(\eta, \sigma ) \cap \{ x=a \} \cap \{|f| = \eta \}.$ Let  
       $$\Psi^{(a)}:  M^{(a)}(\sigma) :  \to \{ a \}\times \Bbb B^2_{\theta } \times \Bbb S^1_{\eta}$$ 
    be    the restriction of $\Psi$ on   $M^{(a)}(\sigma)$.  The intersection 
     $$F_{a,t}=M_t(\sigma)\cap M^{(a)}(\sigma),$$ 
 is nothing but the Milnor fiber  of  the  map germ $f_a : ( \{ a \} \times \Bbb C^2,( a, 0,0)) \rightarrow (\Bbb C,0)$ defined by $f_a(y,z)=f(a,y,z)$ and  $M^{(a)}(\sigma)$ is its  Milnor tube. In Subsection  \ref{Lecar}, we explain  how   the carrousel of D.T. L\^e   is related to  the proof of Theorem \ref{main} .
 }
    \begin{rien} \label{Lecar} {\bf  The carrousel of  D.T. L\^e}
   
 {With the  choice of  coordinates of Subsection 2.1,   $y$ is not a factor of  $f_{a}$. As defined by D.T. L\^e and B. Teissier  (for example see  \cite {Le} or  \cite{T}),  the curve 
    $$\Gamma^{(a)}  :=  H\cap  (\{ a \} \times \Bbb C^2)$$ 
     is  the {\it polar curve} of  $f_{a}$ for the direction $y$  and its  {\it discriminant curve} is:
    $$\Delta^{(a)}  := \Psi(\Gamma^{(a)})=   H' \cap (\{a\} \times \Bbb B^2_{\theta}\times \Bbb B^2_{\eta}).$$
    By definition, $\Delta^{(a)} $ is a plane curve germ at $(a,0,0)$ in  $ \{ a \} \times \Bbb C^2$.
    
       In \cite{Le-2} and \cite {Le},  D.T. L\^e  introduced his  carrousel  construction in order to construct a geometric monodromy  for   hypersurface germs with the help of Puiseux series expansions  of discriminants.    The carrousel  construction,  applied to the germ $f_{a}$, uses a Puiseux series  expansion of  $\Delta^{(a)} $ to study the intersection points  $\Delta^{(a)} \cap {\cal D}_{a,t}$, where  ${\cal D}_{a,t}$ is the disc
  $${\cal D}_{a,t}:=\{a\} \times \Bbb B^2_{\theta} \times \{t\}.$$ 
     It gives   a  decomposition of the solid  torus  $\{ a \}\times \Bbb B^2_{\theta } \times \Bbb S^1_{\eta}$  as a union of Seifert manifolds such that  the iterated torus  link  $ \Delta^{(a)} \cap  (\{ a \}\times \Bbb B^2_{\theta } \times \Bbb S^1_{\eta})$ is a disjoint union of Seifert fibers.    As  $\Psi ^{(a)} $   is a finite  cover   which is ramified over this link,  the pull-back by $\Psi ^{(a)}$ of this Seifert   decomposition of $\{ a \}\times \Bbb B^2_{\theta } \times \Bbb S^1_{\eta},$  provides  a graph manifold  structure on the Milnor tube $ M^{(a)}(\sigma)$. The first return map along the Seifert fibers of the Seifert components gives a quasi-finite monodromy  (for $f_{a}$) on  the Milnor fiber $F_{a,t}.$
   }
\end{rien}

 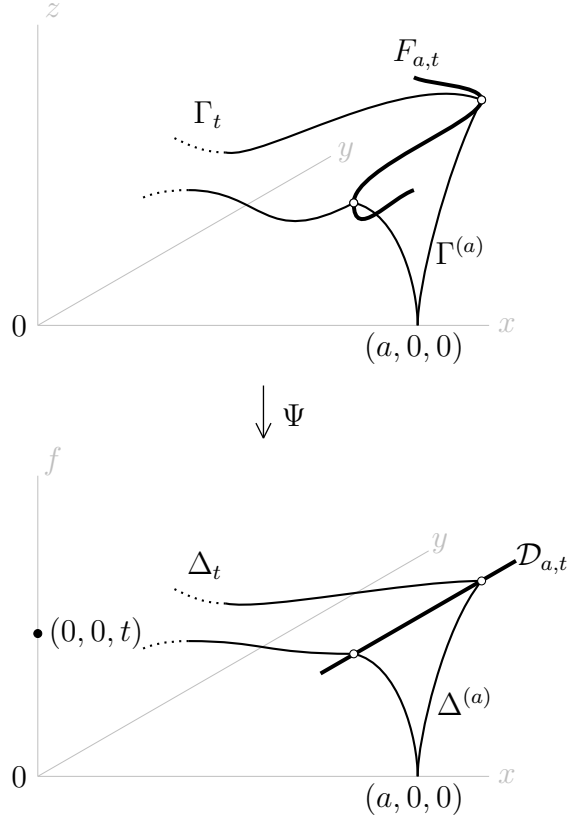
\begin{figure}[ht]
\centering
\begin{tikzpicture} 
\begin{scope}[yshift=6cm]

%EN HAUT

\draw[line width=0.3pt, lightgray]   (0,0)--(6,0);
    \draw[line width=0.3pt, lightgray]   (0,0)--(0,4);
     \draw[line width=0.3pt, lightgray]   (0:0)--(30:4.5);
     
      %\draw[line width=1.5pt,  ]   (20:4)--+(30:3);
      
       \draw[line width=1.5pt,  ]  (5,3.3).. controls (5.2,3.2) and (5.9,3.2) .. (5.9,3);
       \draw[line width=1.5pt,  ]  (4.2,1.63).. controls (4.2,2.1) and (5.9,2.7) .. (5.9,3);
         \draw[line width=1.5pt,  ]  (4.2,1.63).. controls (4.2,1.1) and (4.7,1.7) .. (5,1.8);

    %SUNS
\draw[fill=white] (4.2,1.63)circle(1.5pt);
\draw[fill=white] (5.9,3)circle(1.5pt);
      
%cusp 
 \draw[line width=0.8pt] 
(4.2,1.63) .. controls (4.8,1.5) and (5.05,0.5)  .. (5.05,0) ;
 \draw[line width=0.8pt] 
(5.9,3) .. controls (5.4,2) and (5.05,0.5)  .. (5.05,0) ;
 
 %GAMMA_t
   \draw[line width=0.8pt] 
(4.2,1.63) .. controls (3,1) and (3,1.8)  .. (2,1.8) ;
 \draw[line width=0.8pt, dotted] 
(1.4,1.7) .. controls (1.5,1.75) and (1.7,1.8)  .. (2,1.8) ;

   \draw[line width=0.8pt] 
(2.5,2.3) .. controls (3,2.2) and (4.9,3.4)  .. (5.9,3) ;
 \draw[line width=0.8pt, dotted] 
(2.5,2.3) .. controls (2.2,2.3) and (2,2.4)  .. (1.8,2.5) ;

  %SUNS
\draw[fill=white] (4.2,1.63)circle(1.5pt);
\draw[fill=white] (5.9,3)circle(1.5pt);

%NODES
\node at(2.6,2.8)[left]{$\Gamma_t$};
\node at(6.1,1)[left]{$\Gamma^{(a)}$};
 \node at(5.5,3.6)[left]{$F_{a,t}$};
          
\node at(6.5,0)[left, lightgray]{$x$};
\node at(4.35,2.3)[left, lightgray]{$y$};
\node at(0.2,4.2)[ lightgray]{$z$};
 \node at(0,0)[ left]{$0$};

 \node at(5,-0.3)[]{$(a,0,0)$};

\end{scope}
 
  \node at(3.1,4.8)[right]{$\Psi$};
  \draw[line width=0.5pt]   (3,5.2)--(3,4.5);
   \draw[line width=0.5pt]   (3.1,4.7)--(3,4.5);
     \draw[line width=0.5pt]   (2.9,4.7)--(3,4.5);

%EN BAS

   \draw[line width=0.3pt, lightgray]   (0,0)--(6,0);
    \draw[line width=0.3pt, lightgray]   (0,0)--(0,4);
     \draw[line width=0.3pt, lightgray]   (0:0)--(30:6);
     
      \draw[line width=1.5pt,  ]   (20:4)--+(30:3);
    %  \draw[line width=0.3pt,  ]   (20:4)+(30:1.5)--+(-65:3);
      
%cusp 
 \draw[line width=0.8pt] 
(4.2,1.63) .. controls (4.8,1.5) and (5.05,0.5)  .. (5.05,0) ;
 \draw[line width=0.8pt] 
(5.9,2.6) .. controls (5.4,2) and (5.05,0.5)  .. (5.05,0) ;
 
 %Delta_t
   \draw[line width=0.8pt] 
(4.2,1.63) .. controls (3,1.63) and (3,1.8)  .. (2,1.8) ;
 \draw[line width=0.8pt, dotted] 
(1.4,1.7) .. controls (1.5,1.75) and (1.7,1.8)  .. (2,1.8) ;

   \draw[line width=0.8pt] 
(2.5,2.3) .. controls (3,2.2) and (4.9,2.6)  .. (5.9,2.6) ;
 \draw[line width=0.8pt, dotted] 
(2.5,2.3) .. controls (2.2,2.3) and (2,2.4)  .. (1.8,2.5) ;

\draw[fill ] (0,1.9)circle(1.5pt);

%SUNS
\draw[fill=white] (4.2,1.63)circle(1.5pt);
\draw[fill=white] (5.9,2.6)circle(1.5pt);

%NODES
\node at(2.6,2.8)[left]{$\Delta_t$};
\node at(6.2,1)[left]{$\Delta^{(a)}$};
 \node at(7.2,2.9)[left]{$\cal D_{a,t}$};
          
\node at(6.5,0)[left, lightgray]{$x$};
\node at(5.6,3.1)[left, lightgray]{$y$};
\node at(0.2,4.2)[ lightgray]{$f$};
 \node at(0,0)[ left]{$0$};
                     
\node at(0,1.9)[right]{$(0,0,t)$};
 \node at(5,-0.3)[]{$(a,0,0)$};

 \end{tikzpicture} 
  \caption{The Milnor fiber $F_{a,t}$ as a ramified   cover  over $\cal D_{a,t}$ }
  \label{fig:proof}
\end{figure}

{As already mentioned, the Milnor fiber of the plane curve germ  $f_{a}$ is also included  in the vanishing zone $M_t(\sigma)$ as we have $$F_{a,t}=M_t(\sigma)\cap M^{(a)}(\sigma).$$
The restriction of $\Psi$ to $F_{a,t}$ is a   finite cover over the disc $${\cal D}_{a,t}=\{a\} \times \Bbb B^2_{\theta} \times \{t\}$$ which is ramified over the set of points $H' \cap {\cal D}_{a,t}.$ We obviously have:
  $$H' \cap {\cal D}_{a,t} =  \Delta^{(a)} \cap {\cal D}_{a,t} = \Delta_t \cap {\cal D}_{a,t}.$$
  For this reason, the disc ${\cal D}_{a,t}$ and the restriction $\Psi_{\mid F_{a,t}} \colon F_{a,t} \to {\cal D}_{a,t}$ will play a central role in our proof. This  is also why we need to consider the morphism $\Psi.$ 
  
  Figure \ref{fig:proof} represents  schematically the Milnor fiber $F_{a,t}$ as a ramified  finite cover over the disc ${\cal D}_{a,t}$ and the relative positions of the curves $\Gamma_t$ and $\Gamma^{(a)}$ (resp. $\Delta_t$ and $\Delta^{(a)}$).  

  Here, we have to study the finite cover   $\Psi_t \colon M_t(\sigma) \to \Bbb S^1_{\alpha } \times \Bbb B^2_{\theta} \times \{ t\} , $ which is  ramified over  $H'\cap ( \Bbb S^1_{\alpha } \times \Bbb B^2_{\theta} \times \{ t\})$. Theorem   \ref{par} provides a parametrization  of $H'$ which allows us to obtain Puiseux series   expansions for the following  family of discriminant curves: 
 $$\Delta^{(x)} =   H' \cap (\{x\} \times \Bbb B^2_{\theta}\times \Bbb B^2_{\eta}), \ x\in \Bbb S_{\alpha}^1.$$
  
For us,  {\it  the carrousel in family} is a study  of  the intersections $$  \Delta^{(x)} \cap {\cal D}_{x,t} = \Delta_t \cap {\cal D}_{x,t}, \ x\in \Bbb S^1,$$  with the help of  Puiseux series  expansions  with coefficients in $\Bbb C \{\{x^{1/d}\}\}$. Then, we   trap the braid $H' \cap ( \Bbb S^1_{\alpha } \times \Bbb B^2_{\theta} \times \{ t\})$ in a suitable tubular neighborhood of  an   iterated torus link.   We use   L\^e's ``Swing Lemma"  \cite[Lemma 2.4.7]{L-M-W} (first introduced by D.T. L\^e and B. Perron in \cite {Le-P})   to prove that  that the pull-back (by $\Psi_t$) of  this tubular  neighborhood is a disjoint union of solid tori.}
 
\eject
  \begin{rien} \label{abstract} {\bf  Abstract of the proof of Theorem \ref{main}}
   \vskip.1in
\noindent
   We use the notations of Subsection 4.1.  The proof is  organized as follows. 
\vskip.1in
      \noindent 
 $\bullet$ In Subsection \ref{parametrization}, we describe the braid $H'_t= H' \cap   (\Bbb S^1_{\alpha } \times \Bbb B^2_{\theta} \times \{t\})$ which is the set of ramification values of the finite cover   $\Psi_t \colon M_t(\sigma) \to \Bbb S^1_{\alpha } \times \Bbb B^2_{\theta} \times \{ t\} .$  For this, we consider a  sheet $G$ of  $H$ along the circle $ \Bbb S^1_{\alpha} \times \{0\}\times \{0\}$ as defined in Section 3 and we set $G'=\Psi(G)$.  Using Theorem \ref{par}, we prove the following parametrization result  (Lemma \ref{lem2}): if  $ (x,y,t) \in G'$,  then $y$ satisfies the equality 
   $$y = t^{q/p}\bigg( b\  w(x^{1/d}) \  x^{e /{d'} }     + \sum_{m=1}^{\infty} b_{m}(x^{1/{n'}}) t^{m/pp'} \bigg),\hskip0,5cm \ (\ast)$$ 
where $b \in  \Bbb C^*$, $e \in  \Bbb Z$ and   $d,\ d',\ p,\ q,\ n'$ are positive integers with  $n'=d pp'$ for some $p'$, $w(x^{1/d})= 1+ \sum_{m=1}^{\infty} w_m x^{m/d} \in \Bbb C \{x^{1/d}\} $ and  $ b_{m}(x^{1/n'}) \in    \Bbb C \{\{  x^{1/{n'}}\}\}$. 

We will work inside the solid torus $$T=\Bbb S^1_{\alpha} \times \Bbb B^2_{\theta} \times \{t\}.$$ 
The equality $(\ast)$ implies that $G'_t=\Psi_t(G\cap M_t(\sigma ))$ is a braid in  $T$. We approximate it by the torus link $$App(G'_t)= \{  (x,   b\   x^{e /{d'} }   t^{q/p} , t) ; \ x\in \Bbb S^1_\alpha  \}.$$

 \begin{definition*}   We say that $G'_t $ is {\it the braid of } $G'=  \Psi (G)$,  that $ App(G'_t)$ is {\it the torus link associated to}   $G'$ and that   the pair $( q/p,\ e/d')$ is {\it the pair of  the first exponents  of } $G'$.
  \end{definition*}

      \vskip.1in
      \noindent
$\bullet$ In subsection \ref{pol}, we use the above parametrization to construct a polar decomposition of the solid torus $T$ based on the pairs of  the first exponents  of {all the $G'$'s}.  We  index  the  pairs of first exponents $(q/p,\ e/d')$   {with lexicographic order by decreasing size in $q/p$}. For each of them,  $(q_i/{p_i},\ e_{i,j} /d'_{i,j} )$ where  $ \ 1\leq i \leq k,$ and $\ 1\leq  j \leq l_i $, we construct  a {\it vertical polar zone} $Z_{(i,j)}$ (Definition \ref{polarzone})   such that  $G'_t$  is included in the interior of  $Z_{(i,j)}$ if and only if $ (q_i/{p_i},\ e_{i,j} /d'_{i,j} )$ is the pair of  the  first exponents of $G'$ (Lemma \ref{lem3}). Moreover, $Z_{(1,1)}$  is a solid torus and for all $(i,j)$ not equal to $(1,1)$, the  $Z_{(i,j)}$'s  are concentric thickened tori   which  glue along their boundary components to recover the solid torus $T$.

    \vskip.1in\noindent
 $\bullet$  In subsection \ref{approximation},    we  then define inside $T$ some tubular neighborhoods ${\cal N}(G'_t)$  of   the link  $App(G'_t)$  for all the sheets $G$ of $H$ such that: 
\begin{itemize}
\item[$\ast$] $G'_t \subset {\cal N}(G'_t)$
\item[$\ast$] If $G$ has its pair of first exponents indexed by $(i,j)$, then ${\cal N}(G'_t)$ is included in the interior of $Z_{(i,j)}$.
\item[$\ast$] Let $ \tilde G$ be another sheet of $H$. If   $App(G'_t)=App({ \tilde G}'_t) $, then ${\cal N}(G'_t)= {\cal N}({ \tilde G}'_t)$. Otherwise, ${\cal N}(G'_t)$ and $ {\cal N}({ \tilde G}'_t)$ are disjoint solid tori in $T$ (Lemma \ref{lem4}).
\end{itemize}
 \vskip.1in
 We call  the  solid tori ${\cal N}(G'_t)$  the {\it approximation tori}.      \vskip.1in

\noindent
{\bf  Notation.} Let ${\cal N}(i,j)$ be the union of all the approximation tori  of the sheets which have their first exponents indexed by $(i,j)$.

 By construction the closure of $Z_{(i,j)} \setminus  {\cal  N}_{(i,j)} $ does not meet the set of ramification values  $H'_t$ of $\Psi _t$ and is saturated by $(e_{i,j}, d'_{i,j})$ torus links. The case $e_{i,j}=0$ is not excluded, but  we always have $ 0< \ d'_{i,j} $. It induces a  Seifert  fibration on the closure of  $\Psi ^{-1}_t(Z_{(i,j)} \setminus  {\cal  N}_{(i,j)} )$.

   \vskip.1in\noindent
 $\bullet$  In the last step (Subsections \ref{carrousel in family} and \ref{vertical}), we show (Lemma \ref{lem6}) that $\Psi ^{-1}_t ({\cal N}_{(i,j)})$  is a disjoint  union of solid tori.  The proof of Lemma \ref{lem6} is obtained by applying  L\^e's Swing Lemma   to the family of carrousels for the germs $f_a$ in the direction $y$. Then we can  extend  the Seifert fibration on all the $\Psi ^{-1}_t(Z_{(i,j)} )$. Moreover, we explain in Remark \ref{rk} why the so constructed Seifert manifolds  have oriented basis. It ends the proof of Theorem \ref{main}.  
 
\end{rien} 
  This ends the abstract of the proof of Theorem \ref{main}. 

\begin{rien}  \label{parametrization}
   {\bf  Parametrization of the sheets of  ${\mathbf \Psi (H)}$} \vskip.1in
Let us  recall that $\Psi : (\Bbb C^3,0) \rightarrow (\Bbb C,0)$ denotes the germ defined by $\Psi(x,y,z) = (x,y,f(x,y,z))$.  The critical locus of $\Psi $ is $H = \{ {{\partial f}\over {\partial z}} =0  \} $,  and its discriminant locus is the image $ H' = \Psi (H)$. 
  
Let $G$ be the closure  (in $\Bbb C^3$)  of  a connected component of 
$$ (H \setminus ( \Bbb S^1_{\alpha} \times \{0\}\times \{0\}))  \cap ( \Bbb S^1_{\alpha} \times  \Bbb B^2_{\theta} \times  \Bbb B^2_{\gamma}),$$
{\it i.e.,} $G$ is a sheet of $H$ along the circle $ \Bbb S^1_{\alpha} \times \{0\}\times \{0\}$ as defined in Section 3. We set $G'=\Psi(G)$ and we call $G'$ a {\it sheet} of $ H'=\Psi(H)$ along $ \Bbb S^1_{\alpha} \times \{0\} \times \{0\} .$ 
  \end{rien}
 
\begin{lemma} \label{lem2}
  
 There exist: 
 
 \begin{itemize}
\item[-]  $d,\ n,\ p,\ p',\ q  \in \Bbb N^*$,   where $p$ is prime to $q$  and   $ pp'=n$, 
\item[-]    $e \in \Bbb Z $ and $d' \in \Bbb N^*$ is prime to $e$ (if $e=0$, then $d'=1$),
\item[-] {$ b_{m}(x^{1/n'})  \in   \Bbb C \{\{  x^{1/{n'}}\}\}$, where $ n'=dn$, 
\item[-] $w(x^{1/d})= 1+ \sum_{m=1}^{\infty} w_m x^{m/d}  \in \Bbb C \{x^{1/d}\} $ and $b \in \Bbb C^*$.}
 \end{itemize}
\vskip0,3cm
 such that, if  $ (x,y,t) \in G'$,  then $y$ satisfies the following equality:
  $$y = t^{q/p}\bigg( b\  w(x^{1/d}) \  x^{e /{d'} }     + \sum_{m=1}^{\infty} b_{m}(x^{1/{n'}}) t^{m/pp'} \bigg).   \hskip1cm (\ast)$$ 
  
 \end{lemma}  
\noindent
{{\bf Notation:} From now on, we set $r_{m}:=q/p+m/pp'$. Then,
 $$y =  b\  w(x^{1/d}) \  x^{e /{d'} }   t^{q/p}  + \sum_{m=1}^{\infty} b_{m}(x^{1/{n'}}) t^{r_{m}}.$$}

 \begin{remark*} The integer $d $ is provided by  Theorem \ref{par}. For each  sheet  $G$ of $H$ there exists  such a  minimal $d$ which  depends on  $G$.  Here, for convenience, we will choose a (perhaps greater)  $d$ common to all the sheets of   $H$.
 
 \end{remark*}

\begin{proof}  Let $h:(\Bbb C^3,0)\to(\Bbb C,0)$ be the germ ${{\partial f}\over {\partial z}}$ reduced. Applying Theorem \ref{par} to $h$ provides  $b(x^{1/d},u) \in \Bbb C\{x^{1/d}\}\{u\}$ and $c(x^{1/d},u) \in \Bbb C\{x^{1/d}\}\{u\}$     with $b(x^{1/d},0) \neq 0$ and $c(x^{1/d},0)\neq 0 $, such that we have a parametrization 
 $$   \Bbb S^1_{\alpha^{1/d}}  \times \Bbb B^2_{\epsilon } \rightarrow  G $$ 
 given by 
 $$ (s,u)   \longmapsto (s^d,u^i b(s,u), u^{j} c(s,u)).$$
 Then we obtain $n>1$,  and  $c'(x^{1/d},u) \in \Bbb C\{x^{1/d}\}\{u\}$    with  $c'(x^{1/d},0)\neq 0 $, such that 
 $G'=\Psi(G)$ admits a parametrization   of the form 

 $$(s,u)  \mapsto  (s^d,u^{i} b(s,u), u^n c'(s,u)). \hskip1cm (\ast\ast)$$
 {(The plane curve germ $(\{ f=0\}  \cap (\{ x\}\times \Bbb C^2), (x,0,0))$ is  singular and  $n$ is greater or equal to its  multiplicity.)}
 
 If necessary, we can perform the modification $u= s^{l'}u'$, $l' \in \Bbb N$, to obtain $ l \in \Bbb N$ and  $c_m(s) \in \Bbb C\{ s \}  $  with  $c_0 (0)  \in \Bbb C^{\ast } $,  such that:
  $$t:=u^n c'(s,u)=u'^n  s^l \ c_0(s)\bigg(1+ \sum_{m=1}^{\infty} c_m(s) u'^m\bigg).$$
   There then exist $r(x^{1/d} , u' )\in  \Bbb C\{x^{1/d}\}\{u'\}$ with $ r(0,0)=1 $ and $r_0(x^{1/d}) \in \Bbb C \{ x^{1/d}\} $ with $ (r_0(0))^n=c_0(0)$,   such that 
   $$ t=u'^ns^l(r_0(s))^n(r(s,u'))^n .$$
   We perform the following change of coordinates:$$u_1=u' \  r_0(s)r(s,u'),$$
   and   $(\ast\ast)$ becomes:
    $$(s,u_1)  \mapsto  (s^d,u_1^{i} b'(s,u_1), u_1^n s^l ), $$
      where  $b'(x^{1/d},u_1) \in \Bbb C\{x^{1/d}\}\{u_1\}.$
  
  Now $u_1=s^{-l/n} \ t^{1/n}$ and $ (x ,y,t) \in G'$ satisfies:
  $$y=(x^{-il/nd} \ t^{i/n}) b'(x^{1/d},x^{-l/nd} \ t^{1/n}). \hskip1cm (\ast\ast\ast)$$

  As $x \in \Bbb S^1_{\alpha}$ and $t\in \Bbb B^2_{\eta }$ with $0<\eta << \alpha $, there is no problem of {convergence.  As  $u= s^{l'}u'$ with  $l' \in \Bbb N$  and $u_1=u' \  r_0(s)r(s,u')$ where $ r_0(0)r(0,0)=r_0(0)\neq 0$,  there exists $r'(s)\in  \Bbb C\{ s\},\ r'(s)\neq 0, $ such that  $b'(s,0)=r'(s)b(s,0)$. As  $b(x^{1/d},0) \neq 0$,} we obtain $b\in \Bbb C^*$, $ k\in \Bbb N$ and 
  $w(x^{1/d})= 1+ \sum_{m=1}^{\infty} w_m x^{m/d}  \in \Bbb C \{x^{1/d}\} ,$
    such that:
    $$b'(x^{1/d},0)= b \ x^{k/d} (1+ \sum_{m=1}^{\infty} w_m x^{m/d}) . $$
  
  If we  take  $p$ and $q$ prime to each other such that $q/p=i/n=qp'/pp' $,  $n'=nd$,  $e$ and $d'$ prime to each other such that $e/d'=(n k-i  l)/(nd) $, and  if we write  $(\ast\ast\ast)$  in terms of  increasing powers of $t$,  we obtain  $(\ast)$ of Lemma \ref{lem2},  {\it i.e.,} 
   $$y = b\  w(x^{1/d}) \   x^{e /{d'} }   t^{q/p} + \sum_{m=1}^{\infty} b_{m}(x^{1/{n'}}) t^{r_{m}} .$$ 
   
 This ends the proof of  Lemma \ref{lem2}.
  \end{proof}

\begin{rien} \label{pol} {\bf  The polar decomposition }
 \vskip.1in
  Let us consider the ordered set 
$$\cal P = \bigg\{ \frac{{q_k}}{{p_k}} <   \ldots <  \frac{{q_2}}{{p_2}} <  \frac{{q_1}}{{p_1}}\bigg\}$$  
  of rational numbers $\frac{q}{p}$  such that there exists a sheet $G'$ of  $\Psi(H)$ which admits,  with the notations of \ref{lem2},  a parametrization  of the form:
  
   $$ y=  b\  w(x^{1/d}) \  x^{e /{d'} }   t^{q/p} + \sum_{m=1}^{\infty} b_{m}(x^{1/{n'}}) t^{r_{m}},  $$
   with $ x\in  \Bbb S^1_{\alpha }$ and $t\in \Bbb B^2_{\eta}$.
   
  We denote by $G'_i$ the union of the sheets of $\Psi(H)$ corresponding to the quotient $q_i/p_i$. For  each $i \in \{1,\ldots,k\}$, let 
$$\cal P_i=\{\frac{e_{i,1}}{d'_{i,1}} <  \ldots < \frac{e_{i,j}}{d'_{i,j}} < \ldots < \frac{e_{i,l_{i}}}{d'_{i,l_{i}}}\}$$
\noindent
be the ordered set of rational numbers such that there exists a sheet of $G'_i$   which admits a parametrization of the form: 

   \begin{equation}\label{param}
y=   b\  w(x^{1/d})  \   x^{e_{i,j} /{d'_{i,j}} }   t^{q_i/p_i} + \sum_{m=1}^{\infty} b_{m}(x^{1/{n'}}) t^{r_{m}},
\end{equation}
with $   x\in \Bbb S^1_{\alpha }$ and  $ t\in \Bbb B^2_{\eta}$.

We denote by $G'_{i,j}$ the union of such sheets of $G'_i$. 
          
Let us fix $ a \in \Bbb S^1_{\alpha }$. We consider the  map germ $f_a(y,z)= f (a,y,z)$.
By definition the  set $\cal P$ is the set of polar quotients of $f_a$ for the direction $y$ (for example see \cite {Le}).  We will  follow the classical construction of  \cite {L-M-W}  which furnishes   a  decomposition of the solid torus $T_a = \{a\} \times  \Bbb B^2_{\theta}\times   \Bbb S^1_{\eta}$  into polar zones in bijection with the polar quotients $q_i/p_i$. This decomposition lifts by $\Psi^{(a)}$ to a   decomposition of the exterior  of the link of $f_a $ as a union of Seifert manifolds. But as explained in the abstract of the proof, we have to  define our polar zones $Z_i$ in the solid torus $T = \Bbb S^1_{\alpha} \times  \Bbb B^2_{\theta} \times   \{ t \}$.  The key idea is that the two constructions coincide on the disc ${{\cal D}_{a,t} }= T \cap  T_a$ where they give  a  polar decomposition of ${\cal D}_{a,t}$ as a union of concentric annuli. 

Let us now  define this decomposition of $T$ as  a union of $Z_i$.
 
  For each $i \in \{1,\ldots,k-1\}$, let us choose $s_i \in \Bbb Q$ such that 
 $$ \frac{q_{i+1}}{p_{i+1}}   <  s_i <  \frac{{q_i}}{{p_i}},$$
 \end{rien}
\begin{definition}   The {\it first  polar zone} is the solid torus
 $$Z_1= \{ (x,y,t) \in T / \   |y| \leq \eta^{s_1}\} ,$$
 and  $C(1)= Z_1 \cap {\cal D}_{a,t}$ is the {\it first polar disc}.

 If  $i \in \{2,\ldots, k\}$, the   {\it polar zone} $Z_i$ is the thickened torus defined by:
  $$Z_i = \{ (x,y,t) \in T \  / \  \eta^{s_{i-1}}\leq|y| \leq \eta^{s_i}\} {\hbox{ if } 2\leq i \leq k-1} $$
  {and  
  $$Z_k=\{ (x,y,t) \in T \  / \  \eta^{s_{k-1}}\leq|y| \leq \theta \}.$$}

  The intersection  $C(i)=  Z_i\cap {\cal D}_{a,t}$ is the  associated {\it polar annulus}.
 
 \end{definition}
  
 In $T$, the value of $t\in \Bbb S^1_{\eta}$ is fixed. If $G$ is a sheet of $H$ with  first exponents  $(q_i/p_i, e_{i,j}/d'_{i,j})$, then the braid $\Psi_t (G) = G'_t $  admits a parametrization of the form (\ref{param}) in \ref{pol}.    
 
    To   take account of the first exponent of $x$, we  will refine the polar decomposition of $T$.
 For each $j \in \{1,\ldots,l_i -1 \}$, let us choose a rational number $\nu_{i,j} $ such that 

 $$ \frac{e_{i,j+1}}{d'_{i,j+1}}   <  \nu_{i,j} <  \frac{{e_{i,j}}}{{d'_{i,j}}},$$
 
 There exists $\eta$  sufficiently small, $0< \eta << \theta << \alpha,$  such  that the following   inequalities hold: 
\begin{eqnarray*}
 0 < \eta^{q_1/p_1}\alpha^{\nu_{1,1}} < \eta^{q_1/p_1}\alpha^{\nu_{1,2}} < \ldots < 
 \eta^{q_1/p_1}\alpha^{\nu_{1,l_1 -1}}< \eta^{s_1},
\end{eqnarray*}
for each $ i \in \{2,\ldots,k-1\}$, 
\begin{eqnarray*}
  \eta^{s_{i-1}}<   \eta^{q_i/p_i}\alpha^{\nu_{i,1}} \ldots <  \eta^{q_i/p_i}\alpha^{\nu_{i,l_i -1}} < \eta^{s_i},
\end{eqnarray*}
and
\begin{eqnarray*}
 \eta^{s_{k-1}}<   \eta^{q_k/p_k}\alpha^{\nu_{k,1}} \ldots <  \eta^{q_k/p_k}\alpha^{\nu_{k,l_k -1}} <\theta.
\end{eqnarray*}

 \vskip.1in
\begin{definition} \label{polarzone}
 The {\it vertical polar zones} $Z_{(i,j)}, 1 \leq i \leq k, \  1\leq j \leq l_i,$ are defined as follows: 
 
 \begin{itemize}
\item[$\bullet$] $Z_{(1,1)}$ is the solid torus 
 $$Z_{(1,1)}= \{ (x,y,t) \in T / \   |y| \leq  \eta^{q_1/p_1}\alpha^{\nu_{1,1}} \} ,$$
 \item[$\bullet$]  For $(i,j)$ not equal to $(1,1)$,  $Z_{(i,j)}$ is a thickened torus:
 \begin{itemize}
 \item[$\ast$]   If $1<  i \leq k ,$ 
 
 $$ Z_{(i,1)} = \{ (x,y,t) \in T  \ / \ \eta^{s_{i-1}}\leq |y| \leq  \eta^{q_i/p_i}\alpha^{\nu_{i,1}}\},  $$

\item[$\ast$]  if $1\leq i \leq k,\  j = \{ 2,\ldots,l_i-1\},$ 
$$Z_{(i,j)} = \{(x,y,t) \in T  \ /  \     \eta^{q_i/p_i}\alpha^{\nu_{i,j-1}} \leq |y| \leq \eta^{q_i/p_i}\alpha^{\nu_{i,j}} \}, $$
\item[$\ast$]  if $1\leq i <  k,$
$$Z_{(i,l_i)} = \{(x,y,t) \in T  \ /  \     \eta^{q_i/p_i}\alpha^{\nu_{1,l_i-1}} \leq |y| \leq \eta^{s_i}  \} , $$
\item[$\ast$] and
$$Z_{(k,l_k)}=  \{ (x,y,t) \in T /  \     \eta^{q_k/p_k}\alpha^{\nu_{1,l_k-1}} \leq |y| \leq \theta \}. $$
\end{itemize}
\end{itemize}
 \vskip.1in
The associated refined  {\it polar annuli} are:
$$C(i,j)=Z_{(i,j)} \cap {\cal D}_{a,t}$$

\end{definition}

By construction the torus $T$ is equal to the union of the vertical polar zones $Z_{(i,j)}, \ 1\leq i \leq k ,\ 1\leq j \leq l_i.$ The  intersection of two consecutive (for the lexicographic order on the $(i,j)$) vertical polar zones  is a unique torus which is the common connected component of their boundaries. The intersection between non consecutive vertical polar zones is empty. But, the most important property of the vertical polar zones is given by Lemma \ref{lem3}.

\begin{lemma} \label{lem3}  There exist $\alpha $ and  $\eta$  sufficiently small, $0< \eta << \theta << \alpha,$  such  that  a sheet $G'$ of $H'=\Psi (H)$ has  $(q_i/{p_i},\ e_{i,j} /d'_{i,j} )$ as  pair of  first exponents    if and only if the braid $G'_t = \Psi _t (M_t(\sigma) \cap G)$ is included in the interior of $Z_{(i,j)}.$
\end{lemma}
 \vskip.1in

\begin{proof} By definition, $G'$ has a parametrization of the form (\ref{param}) in \ref{pol}:
 
 $$y=   b\  w(x^{1/d})  \   x^{e_{i,j} /{d'_{i,j}} }   t^{q_i/p_i} + \sum_{m=1}^{\infty} b_{m}(x^{1/{n'}}) t^{r_{m}},$$

 Therefore,  $(x,y,t)\in G'_t$    if and only if
$$\vert y\vert=    \alpha ^{e_{i,j} /{d'_{i,j}} }   \eta ^{q_i/p_i} \  \bigg\vert b\  w(x^{1/d}) \    + \  \sum_{m=1}^{\infty} b_{m}(x^{1/{n'}}) t^{r_m-{q_i/p_i}}  (  x^{ -e_{i,j} /{d'_{i,j}} } )\ \bigg\vert.$$

Then, the inequality 

 $$ \nu_{i,j} <  \frac{{e_{i,j}}}{{d'_{i,j}}} < \     \nu_{i,j-1} $$
implies Lemma \ref{lem3} for the zone $Z_{(i,j)}$ where $ 1\leq i \leq k,\  j = \{ 2,\ldots,l_i-1\}$.

 As $    s_i <  \frac{{q_i}}{{p_i}}  <   s_{i-1}$, the computations are similar for the other vertical polar zones.
\end{proof}

\begin{rien} \label{approximation} {\bf The approximation solid tori}
 \vskip.1in

  Let $G$ be a sheet of $H$ such that  $  G'=\Psi (G)$  is parametrized by   
 \begin{equation*}
y = b\  w(x^{1/d}) \  x^{e /{d'} }   t^{q/p} + \sum_{m=1}^{\infty} b_{m}(x^{1/{n'}}) t^{r_{m}} .
\end{equation*}

  We approximate the braid $G'_t=\Psi_t(G\cap M_t(\sigma ))$ by a torus link $App(G'_t)$ as follows: 
  
\begin{definition*} The {\it link $App(G'_t)$   associated to} the braid  $G'_t $ is the torus link in  $T= \Bbb S^1_{\alpha } \times \Bbb B^2_{\theta} \times \{ t \} $ defined by:

$$App(G'_t)= \{  (x,   b\   x^{e /{d'} }   t^{q/p} , t), \  x\in  \Bbb S^1_\alpha  \}.$$
\end{definition*}
 \vskip.1in
Let   $l$ be the l.c.m. of $d'$ and $p$. Let $a\in \Bbb S^1_{\alpha }$, let $s$ and $\tau$ be such that $s^{d'}=a$ and  $\tau ^p=t$. Let us recall that ${\cal D}_{a,t}$ is the disc $\{a\}  \times \Bbb B^2_{\theta} \times \{ t \}$.
 \vskip.1in
\begin{definition*} The {\it suns of $G'_t$ } are the intersection points ${\cal S}(G'_t) =G'_t\cap {{\cal D}_{a,t}} =\{ (a, b \  \xi \ s^{e}\ \tau^q,t),\ \xi^l=1\}$. 
\end{definition*}
 \vskip.1in

Let $\rho = (e/{d'}+1/{2d})$. 

\begin{definition*} We call {\it approximation solid tori} of $G'_t$ the tubular  neighborhood ${\cal N}(G'_t)$ of $App(G'_t)$  defined by:
$$ {\cal N}(G'_t) = \{(x,y,t) \in T  \  such\  that \     0  \leq |\ y - b\   x^{e /{d'} }   t^{q/p}\  | \leq \eta^{q/p}\alpha^{\rho } \} .$$
\end{definition*}
\end{rien}

\begin{lemma} \label{lem4}  There exist $\alpha $ and  $\eta$  sufficiently small, $0< \eta << \theta << \alpha,$  such  that:
\begin{enumerate}
 \item The intersection $ {\cal N}(G'_t)\cap {{\cal D}_{a,t}} $  consists of  $l$ disjoint  discs of radius equal to $\eta^{q/p}\alpha^{\rho }$  which have  the $l$ suns of $G'_t$ as centers.
\item The braid $G'_t$ is included in  $ {\cal N}(G'_t) $.
\item If   $(q_i/{p_i},\ e_{i,j} /d'_{i,j} )$ is the pair of the first exponents of $G'_t$ then  $ {\cal N}(G'_t) \subset int( Z_{(i,j)})$.
\item   Let $ \tilde G$ be another sheet of $H$. If   $App(G'_t)=App({ \tilde G}'_t) $, then ${\cal N}(G'_t)= {\cal N}({ \tilde G}'_t)$. Otherwise, ${\cal N}(G'_t)$ and $ {\cal N}({ \tilde G}'_t)$ are disjoint solid tori in $T$.
\end{enumerate}

\end{lemma}

 \begin{proof} 
 To obtain (1), it is sufficient to prove that if $ \xi \neq 1$, for a sufficiently small $\alpha$ we have:
 $$ 3\  \eta^{q/p}\alpha^{\rho } <\ |(b-\xi b)| \eta^{q/p}\alpha^{e/d' }.$$ 
 But this inequality is equivalent to:
 $$3\  \alpha^{1/(2d)  }<\ |(b-\xi b)| .$$
As $b$ is a given non zero complex number, it is sufficient to choose $\alpha$  sufficiently small  to obtain (1).

 Let $(s^{d'},y,\tau^p)\in G'_t$ , then: 
 $$y = b\  w(s^{d'/d}) \  s^{e}   \tau ^q + \sum_{m=1}^{\infty} b_{m}(s^{d'/{n'}}) \tau^{pr_{m}} .$$
By construction there exists $w_1(s^{d'/d}) \in \Bbb C\{ s^{d'/d}\}$ such that:

 $$w(s^{d'/d}) -1=s^{d'/d}\   w_1(s^{d'/d}) .$$
For  sufficiently small, $\alpha $ and  $\eta$  with $0< \eta << \theta << \alpha,$   we  have:
$$|y-bs^{e}\tau ^q|= |b\  w_1(s^{d'/d}) \  s^{e+(d'/d)}   \tau ^q + \sum_{m=1}^{\infty} b_{m}(s^{d'/{n'}}) \tau^{pr_{m}}|$$  

$$=\eta^{q/p}\alpha^{e/d'+ 1/d} |b\  w_1(s^{d'/d}) + \sum_{m=1}^{\infty} b_{m}(s^{d'/{n'}})  s^{(-e\ -\ d'/d)}\  \tau ^{-q+pr_{m}}|$$
$$<\ \eta^{q/p}\alpha^{e/d'+ 1/2d} .$$
We then get (2).

To get (3), we show that for  sufficiently small $\alpha $ and  $\eta$  with $0< \eta << \theta << \alpha,$ the distance in ${\cal D}_{a,t}$ between the suns of $G'_t$ and the  boundary components of the polar annulus $C(i,j)$ is bigger  than the radius $ \eta^{q/p}\alpha^{\rho } .$

By construction, we have for $1\leq i \leq k$ and $ j = \{ 2,\ldots,l_i-1\}$: 
$$ C(i,j) = \{(a,y,t) \in {{\cal D}_{a,t}}  \ \hbox{ with }   \     \eta^{q_i/p_i}\alpha^{\nu_{i,j-1}} \leq |y| \leq \eta^{q_i/p_i}\alpha^{\nu_{i,j}} \}, $$
where:
 $$    s_i <  \ \frac{{q_i}}{{p_i}}  <\   s_{i-1} \ ,\ and \    \nu_{i,j} <  \frac{{e_{i,j}}}{{d'_{i,j}}} < \     \nu_{i,j-1}. $$
The distance between a sun of $G'_t$ and the interior circle of  $C(i,j)$ is equal to: 

$$ \eta^{q_i/p_i}\alpha^{(e_{i,j}/d'_{i,j})} (|b|- ( \alpha^{(\nu_{i,j-1})-\ (e_{i,j}/d'_{i,j})} ).$$

This distance, for sufficiently small $\alpha$ and $\eta$, $0< \eta  << \alpha,$  is greater than $ \eta^{q/p}\alpha^{\rho_{i,j} } $ because
 the exponent $\rho_{i,j} = e_{i,j}/{d'_{i,j}}+1/{2d}$ corresponding to a sheet with the pair of the  first exponents equal to $(q_i/{p_i},\ e_{i,j} /d'_{i,j} )$,  is greater than $\frac{e_{i,j}}{d'_{i,j}}$. 
But $ \nu_{i,j} <  \frac{{e_{i,j}}}{{d'_{i,j}}} <  \rho_{i,j}$,  and  similar computations prove that the distance between a sun of $G'_t$ and the exterior circle of $C(i,j)$ is bigger than the radius  $ \eta^{q_i/p_i}\alpha^{\rho_{i,j} }$. The other cases $j=1$ and $j=l_i$ are treated similarly.

Then (3) is done.

Let us now prove (4). When $App(G'_t)=App({ \tilde G}'_t) $, then by definition  ${\cal N}(G'_t)= {\cal N}({ \tilde G}'_t)$.
 
If ${ \tilde G}'_t$ does not have the same pair of first exponents as $G'_t$ then ${\cal N} (G'_t)$ and ${\cal N} ({ \tilde G}'_t)$ are included in the interior of distinct vertical polar zones, they do not meet. 

The last case is when $G'_t$ and ${ \tilde G}'_t$ have the same   pair of   first exponents $(q/p,e /d' )$, but distinct associated torus links. If  $(s^{d'},y,\tau^p)\in G'_t$, then: 
 $$y = b\  w(s^{d'/d}) \  s^{e}   \tau ^q + \sum_{m=1}^{\infty} b_{m}(s^{d'/{n'}}) \tau^{pr_{m}} .$$
If    $(s^{d'},y,\tau^p)\in { \tilde G}'_t$ , then: 
 $$y =\tilde{b}\  \tilde{w}(s^{d'/d}) \  s^{e}   \tau ^q + \sum_{m=1}^{\infty} \tilde{b}_{m}(s^{d'/{n'}}) \tau^{pr_{m}} ,$$
where $\tilde{b}  \in \Bbb C^* $ and  $\tilde{b}\neq \xi  b, $  for  all  $\xi $ such that  $ \xi^l=1.$
But the minimal value of $ \{ |\tilde{b} -\xi b|,\ \xi^l=1\}$ is well defined. With  computations similar to those performed to obtain  points (1) and (2), we can choose  sufficiently small $\alpha$ and $\eta$, $0< \eta  << \alpha,$ such that the distances between the suns of $G'_t$ and ${ \tilde G}'_t$ are bigger than $3 \eta^{q/p}\alpha^{\rho }$. This proves that ${\cal N} ({ \tilde G}'_t)$ and ${\cal N} (G'_t)$ are disjoint.

 But the trivial projection  of $T= \Bbb S^1_{\alpha } \times \Bbb B^2_{\theta }\times \{ t \}$ on $ \Bbb S^1_{\alpha} $ restricted on  ${\cal N} (G'_t)$ is a fibration with the discs $ {\cal N}(G'_t)\cap {{\cal D}_{a,t}} $ as fibers. Then the tubular neighborhoods ${\cal N}(G'_t)$ are unions of disjoint solid  tori in $T$.
 
This ends the proof of Lemma \ref{lem4}.
\end{proof}
 
 Lemma \ref{lem4} allows us to define the solar discs.
 
 \begin{definition} Let $s$ and $\tau$ be such that $s^{d'}=a$ and  $\tau ^p=t$. If $G$ is a sheet of $H$ and  $G'=\Psi (G)$, the solar discs  associated to $G$ are the $l$ disjoint discs $ {\cal N}(G'_t)\cap {{\cal D}_{a,t}} $ centered at the suns ${\cal S}(G'_t) =G'_t\cap {{\cal D}_{a,t}} =\{ (a,b \ \xi \ s^{e}\ \tau^q,t),\ \xi^l=1\}$ of $G'_t$.
\end{definition}

\begin{lemma} \label{lem5} Let $D_G$ be a solar disc of $G$, then   $\Psi ^{-1}_t (D_G))$  is a disjoint union of discs.
\end{lemma}

\noindent
  Subsection  \ref{carrousel in family} contains the proof of Lemma  \ref{lem5}.
  
\begin{rien} \label{carrousel in family}
{\bf  Carrousel in family}
 \vskip.1in

{Recall that  $\Psi_t \colon M_t (\sigma ) \to (\Bbb S^1_{\alpha} \times \Bbb B^2_{\theta} \times \{t\} )$ denotes the restriction of $\Psi$ (where $\Psi(x,y,z) = (x,y,f(x,y,z))) $  on $M_t (\sigma)$.  We use again the notations of Subsection \ref{abstract}.  Let us fix $ a \in \Bbb S^1_{\alpha}$.  By  construction,  the Milnor fiber of the map germ $f_a : ( \{ a \} \times \Bbb C^2,( a, 0,0)) \rightarrow (\Bbb C,0)$ defined by $f_a(y,z)=f(a,y,z)$  is  $$F_{a,t} =M_t(\sigma )\cap \{ x=a \}.$$  
  
  Let  $$\psi_a  \colon F_{a,t} \rightarrow {\cal D}_{a,t} .$$ 
  be the restriction of $\Psi $ on $F_{a,t}$.} Then Lemma \ref{lem5} is equivalent to:

 \vskip.1in
  {\bf  Claim.}  Let $D_G$ be a solar disc of $G$, then   $\psi ^{-1}_a (D_G))$  is a disjoint union of discs.
 \vskip.1in

\begin{proof}[Proof of the claim] Let $\delta$ be a irreducible component of the discriminant $\Delta_a$ which is included in $G'=\Psi (G)$. Then a Puiseux series  expansion of $\delta$ is  given by:
 $$y = b\  w(s^{n}) \  s^{ed''}   t^{q/p} + \sum_{m=1}^{\infty} b_{m}(s) t^{r_{m}} .$$
Where  $s$  and $d''$ satisfy the following equalities:  $s^{nd}=a$ and  $d'd''=nd$. Moreover,  the suns of $\delta$ as defined in \cite[Definition 2.4.3]{L-M-W} are the following  $p$ points  of ${\cal D}_{a,t}$: $\{ (a,b\  w(s^{n}) \  s^{ed''}   \tau^{q}, t),\ \tau^p=t \}$. 
 A solar ``polar" disc $D$ is defined in \cite[Definition 2.4.6]{L-M-W}, and  \cite[Lemma 2.4.7]{L-M-W}  states that $\psi_a^{-1}(D) $ is a disjoint union of discs. 
This uses L\^e's swing. 
Our polar disc $D_G$  takes account of the coefficients parametrized by $x$ via $w(x^{1/d})$ and  is slightly different from $D$. But
 we  can consider the  curve  $\delta '$ having  $y = b \  s^{ed''}   t^{q/p} $ as Puiseux series  expansion in $\{ a \} \times \Bbb C^2$. If we  use  the curve $\delta '$ in the proof of \cite[Lemma 2.4.7]{L-M-W}  in place of $\delta _0$, we obtain, with exactly the same arguments, that   $\psi ^{-1}_a (D_G))$  is a disjoint union of discs. 
 This proves the claim. 
 \end{proof}
\end{rien}
 \vskip.1in
\begin{remark} In \cite {C}, C. Caubel proves a very  general version of  L\^e's swing. In particular let $D$ be a subdisc of a polar annuli $C(i,j)$. We say that $D$ is marked if it contains points of {the discriminant curve $\Delta ^{(a)}$ in its interior, but the boundary of $D$ does not meet $\Delta^{(a)}$}. \cite[Proposition 2.4 ]{C} implies that  if  $D$  is  a marked subdisc contained in a sector, in  $C(i,j)$, of angle $\theta$ with $\theta <\ 2\pi (q_i/p_i+1/2p_i)$, then  $D$ can be swung. 

Then, in the case of a germ  defined on $(\Bbb C^2,0)$, as {$f_a$} in  our  case, we obtained  (as proved in  \cite[2.4.12]{L-M-W}), that $\psi ^{-1}_a (D)$ is a disjoint union of discs. By definition  our polar disc $D_G$ is contained in such a sector.

\end{remark}

 \vskip.1in
\begin{rien} \label{vertical} {\bf  Vertical monodromy}
 \vskip.1in
Let $  p$  be the restriction on $M_t(\sigma )$  of the projection on  the $x$-axis:
$$p:  M_t(\sigma )\ \rightarrow \ \Bbb S^1_{\alpha }.$$

In  Subsection \ref{axis1} we choose a generic $x$-axis such that ${\it p}$ is a submersion on $M_t(\sigma )$ when $t \in \Bbb S^1_{\eta},\  0<\eta<<\alpha.$ Then $p$ is a differentiable fibration  with fiber $F_{a,t}$ and $M_t(\sigma )$ is the mapping-torus of  a diffeomorphism $h : F_{a,t} \rightarrow F_{a,t}$. Following the terminology introduced by D. Siersma in \cite{S},  $h$  is a representative of the {\it vertical monodromy} for $\sigma$.

Let ${\cal N}(i,j)$ be the union of all the approximation tori  of the sheets which have their first exponents indexed by $(i,j)$.

\end{rien}

 \vskip.1in
\begin{lemma} \label{lem6} Each   $\Psi ^{-1}_t ({\cal N}_{(i,j)})$  is a disjoint union of solid tori.
\end{lemma}

\begin{proof}[Proof of Lemma \ref {lem6}]  By construction the boundary of  $\Psi ^{-1}_t ({\cal N}_{(i,j)})$ meets $\{ x=a \}$ transversally for all $a\in \Bbb S^1_{\alpha}$. Then,  the restriction  ${ \it p_{i,j}}$  of $ { \it p}$  on  $\Psi ^{-1}_t ({\cal N}_{(i,j)})$ is a fibration. But the fibers of this restriction consist of a disjoint union of $\Psi ^{-1}_t (D_G)$ for all the polar discs $D_G$  of the sheets $G'= \Psi (G)$ having $(q_i/p_i, e_{i,j}/d'_{i,j})$ as pair of first exponents. Lemma \ref {lem5} implies that the fibers of $ { \it p_{i,j}}$  are a disjoint union of discs. Then $\Psi ^{-1}_t ({\cal N}_{(i,j)})$ is the mapping torus of a disjoint union of discs, it is a disjoint union of solid tori.
\end{proof}

Lemma \ref {lem6} is the key lemma which  enables one to conclude. 

By construction the closure of $Z_{(i,j)} \setminus  {\cal  N}_{(i,j)} $ does not meet the ramification value $H'_t$ of $\Psi _t$ and is saturated by $(e_{i,j}, d'_{i,j})$ torus links. It induces a Seifert fibration on the closure of  $\Psi ^{-1}_t(Z_{(i,j)} \setminus   {\cal N}_{(i,j)} )$. Moreover,  the Seifert   fibers are, by construction, transverse to the fibers of ${\it p}$. Then, Lemma \ref {lem6} allows us to extend the Seifert fibration  on the  disjoint union of solid tori  $\Psi ^{-1}_t ({\cal N}_{(i,j)}).$ The connected components of $\Psi ^{-1}_t (D_G))$ being the meridian discs of the tori $\Psi ^{-1}_t ({\cal N}_{(i,j)})$, there is no singular   fiber  in the constructed Seifert fibration on  $\Psi ^{-1}_t(Z_{(i,j)} ) $ and the possible exceptional  Seifert fibers  are the cores of the tori $\Psi ^{-1}_t ({\cal N}_{(i,j)})$  or are in $\Psi ^{-1}_t  (\Bbb S^1_{\alpha } \times \{ 0 \}  \times \{ t  \} )$. The union  along their boundary components of the Seifert manifolds $\Psi ^{-1}_t(Z_{(i,j)})$, for all $(i,j)$,  gives a decomposition  of 
 $M_t (\sigma ) = \Psi ^{-1}_t  ( \Bbb S^1_{\alpha } \times \Bbb B^2_{\theta} \times \{ t \} )$ as a graph manifold.

 \begin{remark} \label{rk} The above constructed Seifert  fibers  of $M_t(\sigma )$ define a quasi-finite vertical monodromy which preserves the orientation of the oriented Milnor fiber $F_{a,t}$. It implies that the obtained Seifert pieces of  $M_t(\sigma )$ have oriented basis.
 \end{remark}
   This ends the proof of   Theorem \ref {main}. \hfill$\Box$

 \section{A topological characterization of isolated singularities } 
In this section, we prove the following topological characterization of isolated singularities, which was the first motivation of this work. 
 
 \vskip.1in
\begin{theorem} \label{maintop} Let $f:({\Bbb C}^3,0) \longrightarrow ({\Bbb C},0)$ be a reduced holomorphic
germ. Unless $f$ is irreducible and $L_t$ is  a lens space,  the following assertions are equivalent.
\begin{itemize}
\item[(i)] $f$ is either smooth or has an isolated singularity at $0$.
\item[(ii)] The boundary $L_t, \ t \not=0$, of the Milnor fiber of $f$ is homeomorphic to the link $\overline{L_0}$ of the
normalization of $f^{-1}(0)$.
\end{itemize}
\end{theorem}

The degenerating case when $f$ is irreducible and $L_t$ is a lens space remains open.

If $f$ is reducible,  $L_t$ is not  homeomorphic to $\overline  L_0$. Indeed,  by definition, the number of connected components of $\overline  L_0$  equals  the number of irreducible components of $f$, but $L_t$ is always a connected manifold (Corollary \ref{connected}). On the other hand, $\overline  L_0$ is an irreducible $3$-dimensional  \cite[Theorem 1]{N}.   Then, it suffices to prove the theorem when $f$ is an irreducible germ and $L_t$ is an irreducible $3$-dimensional manifold. From now on,  we assume that $f$ is irreducible.

Before proving the theorem, we will establish some basic properties of $L_t$.

 \vskip.1in
 
 \begin{proposition} \label{irred}
The  trunk $\overline{N_0}$   and the vanishing zone $M_t$   are irreducible $3$-manifolds. 
 \end{proposition}
 Recall that a $3$-manifold $M$ is irreducible if every embedded $2$-sphere in $M$ is the boundary of a $3$-ball.
  \vskip.1in
\begin{proof}[Proof of Proposition \ref{irred}] It suffices to prove that every connected component $\overline{W}$ of $\overline{N_0}$ is irreducible. Let $(S,p)$ be an irreducible component of $\overline{F_0}$ whose link contains $\overline{W}$, and set  $\gamma = \overline{\Sigma(f)} \cap S$. Then $\overline{W}$ is the complement of a tubular neighborhood of the link of the complex germ of curve $(\gamma,p)$ in the link of the normal complex surface singularity  $(S,p)$. Therefore $\overline{W}$ is irreducible (see \cite[9.2, Corollary J]{M-P-W-2}).
 
 According to Subsection  \ref{vertical}, each connected component $M_t(\sigma)$ of the vanishing zone $M_t$ is fibered  over the circle $\Bbb S^1$ with a connected and orientable fiber  not diffeomorphic to the $2$-sphere. Therefore  $M_t(\sigma)$ is irreducible   (see \cite[9.1, Lemma A]{M-P-W-2}).  
 \end{proof}

  \begin{corollary} \label{irr} Assume  that $f$ irreducible and that $\overline{N_0}$ is not a solid torus. Then $L_t$ is an irreducible $3$-dimensional manifold.
  \end{corollary}
\begin{proof} It is an easy consequence of the following general principle, which is a consequence of \cite{W} : let $(M_i), i=1,\ldots,k$ be a finite collection of Seifert manifolds with non empty boundary, none of them being a solid torus. Let $M$ be constructed by gluing the $M_i$'s along boundary tori. Then $M$ is irreducible.  
\end{proof}

Notice that $\overline{N_0}$ is a solid torus if and only if the minimal resolution graph of $(\overline  L_0, \overline{\Sigma}(f))$ is a bamboo with an arrow at one of its extremities. 
 \vskip.1in
\begin{remark} In fact,    when $f$ is irreducible, $L_t$ is a reducible $3$-dimensional manifold   if and only if  $\overline{N_0}$ is a solid torus and a Seifert fiber on the boundary of $M_t$ is a meridian of $N_t \cong \overline{N_0}$.
\end{remark}
 \vskip.1in
Let $M$ be an  irreducible graph  manifold. We denote by ${\cal T}(M)$ the  separating  family of  the minimal Jaco-Shalen-Johannson decomposition of $M$ and by $\sharp {\cal T}(M)$ the  cardinality of ${\cal T}(M)$.
When $M$ has empty boundary, we  denote by ${\cal G}(M)$ the normalized plumbing graph of $M$ as defined in   \cite{N}.

\begin{proposition} \label{complexity} Assume that the germ $f$ is irreducible and that $L_t$ is an irreducible $3$-dimensional manifold. Then

$$ \sharp {\cal T}(\overline{L_0}) \leq \sharp {\cal T}(\overline{N_0}) \leq  \sharp {\cal T}(L_t) , \hbox{ and }$$
$$rank \ H_1({\cal G}(\overline  L_0), \Bbb Z)  \leq  rank \ H_1({\cal G}(L_t), \Bbb Z) $$

\end{proposition} 

\begin{proof} When $\overline{N_0}$ is a solid torus, then $\overline{L_0}$ is a lens space. Then $rank \ H_1({\cal G} (\overline  L_0), \Bbb Z) =0$,  $\cal T(\overline{N_0}) = \emptyset$ and the two inequalities hold.

Assume that $\overline{N_0} \cong N_t$ is not a solid torus, then 
$$\cal{T}(M_t) \cup \cal{T}(N_t) \subset \cal{T}(L_t) $$

On the other hand, one has $ \sharp {\cal T}(\overline{L_0}) \leq \sharp {\cal T}(\overline{N_0})$ as the closure of $\overline{L_0} \setminus \overline{N_0}$ is a disjoint union of solid tori.

When $M$ is an irreducible graph manifold without boundary, we denote by ${\cal G}'(M)$ the graph ${\cal G}(M)$ without weights, and we extend the definition of ${\cal G}'(M)$ to the case when $M$ has a non empty boundary by symbolizing each boundary component of $M$ by a free edge. 

For example, the graph  ${\cal G}'(\overline{N_0})$ is obtained from   the normalized  plumbing graph with arrows of the pair $(\overline{L_0}, \overline{\Sigma}(f))$ by removing the weights (genus and Euler classes) and by replacing each arrow by a free edge. 

According  to Theorem  \ref{main1}, $L_t$ is obtained from $N_t \cong \overline{N_0}$ and $M_t$ by  gluing together these  two manifolds along their boundary components. Then the graph ${\cal G}'(L_t)$ is obtained from the two graphs ${\cal G}'(\overline{N_0})$ and  ${\cal G}'(M_t)$ by identifying the free edges corresponding to the glued boundary components. This proves the  second inequality.
\end{proof}

\begin{proof}[Proof of Theorem \ref{maintop}]  (i) $\Rightarrow$ (ii) follows from Milnor's theory (\cite{Mi}). 
\par
To prove (ii) $\Rightarrow$ (i), let us assume that $f$ is neither  smooth nor has an isolated singularity at $0$. As mentioned at the beginning of the Section, it suffices to prove  (ii) $\Rightarrow$ (i) when $f$ is an irreducible germ  and $L_t $ is an irreducible $3$-manifold.

If the trunk  $\overline{N_0} \cong N_t$ is a solid torus, then $\overline{L_0}$ is a lens space. But, we have assumed that $L_t$ is not a lens space, so $L_t$ is not homeomorphic to $\overline{L_0}$.
 
Now, assume   that the trunk  $\overline{N_0} \cong N_t$ is not a solid torus. Then $L_t$ is obtained as the union of  the two irreducible manifolds $M_t$ and $N_t$ along their boundaries  (Theorem \ref{main1}), none of them being a solid torus. Therefore $L_t$ is irreducible (Corollary \ref{irr}).

Assume first that there exists a connected component $M_t(\sigma)$ of $M_t$ whose boundary is not connected.   Gluing the manifold $M_t(\sigma)$ to the trunk $N_t$ increases the number of cycles in the normalized  plumbing graph ${\cal G}(\overline  L_0)$. Therefore, 
$$rank \ H_1({\cal G}(\overline  L_0), \Bbb Z)  <  rank \ H_1({\cal G}(L_t), \Bbb Z) $$
  and   $\overline{L_0}$ is not homeomorphic to $L_t$.

We now assume that  each connected component $M_t(\sigma) $  of $M_t$ has a connected boundary,  {\it i.e.,} $M_t(\sigma) \cap N_t$ consist of a single torus.

When   $ \sharp {\cal T}(\overline{L_0}) <   \sharp {\cal T}(L_t) $, then $\overline{L_0}$ is not homeomorphic to $L_t$. Otherwise, the equality $ \sharp {\cal T}(\overline{L_0}) =   \sharp {\cal T}(L_t) $ implies that the connected components of $M_t$ are all Seifert manifolds and that the Seifert structure induced on the boundary components are homological to that of $N_t$. 

We now use the following:

\begin{remark*} Let $M$ be an irreducible  orientable $3$-dimensional manifold whose Jaco-Shalen-Johannson decomposition admits only Seifert pieces with orientable basis. Assume that $M$ is  neither  diffeomorphic  to a lens space nor to a solid torus. Then, according to the classical classification of irreducible $3$-dimensional manifolds (see \cite{Ha}), the following two numbers  are numerical  invariants of the homeomorphism class of $M$: 

\begin{enumerate}
\item  The sum $g(M)$ of the genus of the  bases of the Seifert pieces of $M$ in any Jaco-Shalen-Johannson decomposition of $M$,
\item the global number $s(M)$ of exceptional Seifert  fibers  in the minimal decomposition of $M$.
\end{enumerate}
\end{remark*}
 \vskip.1in
 
 Let $r$  be the number of irreducible components of $\Sigma(f)$. As each connected component $M_t(\sigma)$ of $M_t$ has a connected boundary, then $r$ is also the number of irreducible components of the curve  $\overline{\Sigma}(f)$.  Therefore the trunk $\overline N_0 \cong N_t$ has $r$ boundary components (Corollary \ref{sigma}), and $\overline{L_0}$ is obtained by gluing  $r$ solid tori along the $r$ boundary components of $ \overline N_0$. We then have: 
  $$ g(\overline{L_0}) = g( \overline N_0) \hbox{ and } s(\overline{L_0}) \leq s(\overline N_0) +r \hskip 1cm (\ast)$$

Let $\sigma$ be an irreducible component of $\Sigma(f)$. Let $p : M_t(\sigma) \rightarrow \Bbb S^1_{\alpha}$ be the 
locally trivial fibration with fiber $F_{a,t}$ and monodromy $h : F_{a,t} \rightarrow F_{a,t}$ defined in Subsection \ref{vertical}.

If the transversal section of $F_0$ at a point  of $\sigma \setminus \{0\}$ is the ordinary quadratic germ, then the Milnor fiber $F_{a,t}$ is an annulus $[-1,+1] \times \Bbb S^1$. As $M_t(\sigma)$ has a connected boundary, then $h : [-1,+1] \times \Bbb S^1 \rightarrow [-1,+1] \times \Bbb S^1$ is isotopic to the diffeomorphism $h(t,z)=(-t,\bar{z})$ and its mapping torus $M_t(\sigma)$ is the so-called Seifert Q manifold (\cite{W}), which has two exceptional fibers and base a disk.

In all other cases, $\chi(F_{a,t}) < 0$. Then $M_t(\sigma)$ has either $g(M_t(\sigma))>0$ or at least two exceptional fibers, {\it i.e.,} $s( M_t(\sigma))\geq 2$.

If there exists $\sigma$ such that $g(M_t(\sigma))>0$, then $g(L_t)> g(N_t) = g( \overline N_0) =g(\overline{L_0})  $, so $L_t$ is not homeomorphic to  $\overline{L_0}$. 

Otherwise, each $M_t(\sigma)$ has at least $2$ exceptional fibers, and 
$$s(L_t) \geq s(N_t) + 2r$$
Then $(\ast)$ implies $s(L_t) > s(\overline{L_0)}$ and $L_t$ is not homeomorphic to  $\overline{L_0}$.      
\end{proof}

Theorem \ref{maintop} remains open when $f$ is irreducible and $L_t$ is a lens space. The following proposition shows that, in fact, this case concerns a very special family of singularities. Recall that the $\overline{K_0}$ denotes the link of the curve $\overline{\Sigma}(f) $ in the link $\overline{L_0}$ of the normalization $\overline{F_0}$ of $F_0$.

 \vskip.1in
\begin{proposition} \label{lens} Let $f:({\Bbb C}^3,0) \longrightarrow ({\Bbb C},0)$ be a reduced holomorphic
germ such that $f$ is irreducible and $L_t$ is  a lens space. Then 
\begin{enumerate}
\item The trunk $\overline N_0$ is a solid torus,  $\overline{L_0}$ is a lens space, $\overline{\Sigma}(f) $ is an irreducible  curve germ   and  the minimal  resolution graph of the pair $(\overline{F_0},\overline{\Sigma}(f))$ is a bamboo with an arrow at one of its extremities,  
 \item $M_t$ is connected with a connected boundary.
\end{enumerate}
 \end{proposition}
\begin{proof}[Proof of Proposition \ref{lens}]
Let $\sigma$ be a component of $\Sigma(f)$. According to  \ref{vertical}, $M_t(\sigma)$ is  fibred over the circle with fiber $F_{a,t}$. As $F_{a,t}$  is not a {disc,} then   $M_t(\sigma)$ is not a solid torus. 

Let $T$ be a connected component of $\partial N_t =  \partial M_t$. As the connected components of $M_t$ are  irreducible  manifolds (Proposition \ref{irred}) none of them being  a solid torus, then $T$ is incompressible in $M_t$ (see \cite{M-P-W-2}, 9.1, prop. D). Now, as the trunk  $\overline{N_0} \cong N_t$ is irreducible (\ref{irred}),  if it were not a solid torus,  $T$ would also be  incompressible in $N_t$ (see again \cite{M-P-W-2}, 9.1, prop. D). Then,  van Kampen's theorem and Dehn's   lemma would   imply  that $T$ is incompressible in $L_t$. But  a torus embedded in a lens space is always compressible. Hence  $\overline{N_0} \cong N_t$ is a solid torus and then  the minimal  resolution graph of the pair $(\overline{F_0},\overline{\Sigma}(f))$ is a bamboo with an arrow at one of its extremities. 
  It  follows immediately that $\overline{L_0} $ is a lens space. According to \ref{sigma}, the curve $\overline{\Sigma}(f)$ is irreducible in $\overline{F_0}$. Therefore $\Sigma(f)$ is also irreducible. 
 
 As  the trunk  $\overline{N_0} \cong N_t$ is a solid torus,  the vanishing zone $M_t$ is connected with a connected boundary because $\partial N_t = \partial M_t$.   
\end{proof}

\vskip.1in

\noindent {\bf Adresses.}

\vskip.1in

\noindent Fran\c coise Michel / Laboratoire de Math\' ematiques Emile
Picard  /
Universit\' e Paul Sabatier / 118 route de Narbonne / F-31062
Toulouse / FRANCE

e-mail: fmichel@picard.ups-tlse.fr

\vskip.1in

\noindent Anne Pichon / Aix Marseille Universit\'e, CNRS, Centrale Marseille / I2M, UMR 7373 \ 13453 Marseille \ FRANCE

e-mail: anne.pichon@univ-amu.fr

\vskip.5in

\end{document}